\documentclass[12pt, leqno]{article}
\usepackage{amsmath}
\usepackage{amsfonts}
\usepackage{amssymb}
\usepackage{graphicx}
\usepackage{epsfig}     
\usepackage{color}
\usepackage{color}
\newtheorem{thm}{Theorem}[section]
\newtheorem{cor}[thm]{Corollary}
\newtheorem{lem}{Lemma}[section]
\newtheorem{prop}[thm]{Proposition}
\newtheorem{rem}[thm]{Remark}
\numberwithin{equation}{section}

\newcommand{\dx}{\,{\rm d}x}

\def\supp{\mathrm{supp}} 
\newcommand{\RR}{\mathbb{R}}
\newcommand{\re}{\mathbb{R}}
\newcommand{\ren}{\mathbb{R}^N}

\newcommand{\la}{\lambda}
\newcommand{\ve}{\varepsilon}
\newcommand{\vp}{\varphi}

\newcommand{\dint}{\displaystyle\int}
\newcommand{\diint}{\displaystyle\iint}
\renewcommand{\dfrac}{\displaystyle\frac}
\def\qed{\,\unskip\kern 6pt \penalty 500
\raise -2pt\hbox{\vrule \vbox to8pt{\hrule width 6pt
\vfill\hrule}\vrule}\par}
\definecolor{darkblue}{rgb}{0.05, .05, .65}
\definecolor{darkgreen}{rgb}{0.05, .70, .05}
\definecolor{darkred}{rgb}{0.8,0,0}
\def\qed{\unskip\kern 6pt \penalty 500
\raise -2pt\hbox{\vrule \vbox to8pt{\hrule width 6pt
\vfill\hrule}\vrule}\par}
\begin{document}
\title{\textbf{Regularity of solutions of the \\ fractional  porous medium flow }}
\author{\Large Luis Caffarelli \footnote{caffarel@math.utexas.edu},\\
\Large Fernando Soria \footnote{fernando.soria@uam.es},
 \\\Large Juan Luis Vazquez\footnote{juanluis.vazquez@uam.es}\\}

\date{\today}
\maketitle

\begin{abstract}
We study a porous medium equation with nonlocal diffusion effects
given by an inverse fractional Laplacian operator. More precisely,
$$
u_t=\nabla\cdot(u\nabla (-\Delta)^{-s}u), \quad  \ 0<s<1.
$$
The problem is posed in $\{x\in\ren, t\in \re\}$ with nonnegative initial data $u(x,0)$ that are integrable and decay at infinity. A previous paper has established the existence of mass-preserving, nonnegative weak solutions satisfying energy estimates and finite propagation. Here we establish the boundedness and $C^\alpha$ regularity of such weak solutions

\end{abstract}

\newpage


\section{Introduction}
\label{sec.intro}

This paper is devoted to study the regularity properties of weak solutions of a model of porous medium equation that includes nonlocal effects through an integral relation of pressure to density. This allows to account for long-range effects.

Let us recall the typical derivation of the porous medium equation, cf. \cite{Ar, Vapme}. We consider a  gas propagating in a homogeneous  porous medium; its dynamics is described by first assuming conservation of mass
\begin{equation*}
\partial_t u + \nabla \cdot  ( {\bf v} u)=0\,,
\end{equation*}
where $u(x,t)\ge 0$ denotes the density of the gas and ${\bf v}(x,t)$ is the (locally averaged) velocity. We then postulate that the
motion proceeds according to Darcy's law so that ${\bf v}=-\nabla p$, where the velocity potential is interpreted as a pressure. Finally, some  barotropic state law for gases implies that $p$ is a monotone function of $u$, $p=f(u)$. In this way we get the equation
\begin{equation*}
\partial_t u= \nabla \cdot  (u\nabla f(u)).
\end{equation*}
The simplest case (called isothermal) is $p=u$ and in that case we arrive at $\partial_t u=c\Delta u^2$, which appears in a different context as a model in groundwater infiltration, Boussinesq's equation  \cite{Bear, Bou03} .

The novelty of our present model consists is relating $p$ to $u$ through a linear  integral operator that makes a kind of average of the space distribution $u(\cdot,t)$,
\begin{equation}
p(x,t)={\cal L} u(x,t), \qquad {\cal L} u(x,t):=\int L(x-y)\,u(y,t)\,dy
\end{equation}
More in particular, the positive  kernel $L$ is locally integrable and decays slowly at infinity to represent ``long-range '' interactions. To  be specific, we will work in $\ren$, we will take $L(x)=c|x|^{-N+2s}$, which is equivalent to saying that $p$ is given as an inverse fractional Laplacian, i.\,e., $p=(-\Delta)^{-s} u$, and we consider $0<s<1$.

In a previous paper \cite{CaVa09} we have introduced this model and proved existence of weak solutions for the Cauchy problem
\begin{equation}\label{eq1}
u_t=\nabla\cdot(u\nabla p), \quad p={\cal L}u=(-\Delta)^{-s}u, \ 0<s<1\,,
\end{equation}
posed for $x\in \RR^N$, $N\ge 1$, and $t>0$, with  initial conditions
\begin{equation}\label{eq.ic}
u(x,0)=u_0(x), \quad x\in \RR^N,
\end{equation}
where $u_0$ is a nonnegative and integrable function in $\RR^N$ decaying as $|x|\to\infty$.

Let us point out that equations of the more general form $u_t=\nabla \cdot (\sigma(u)\nabla {\cal L}u)$ have appeared recently in a number of applications in particle physics. Thus, Giacomin and Lebowitz consider in   \cite{GL1} a lattice gas with general short-range interactions
and a Kac potential $ J_\gamma(r)$ of range $\gamma^{-1}$, $\gamma>0$ . Scaling
spacelike with $\gamma^{-1}$ and timelike with $\gamma^{-2}$, and passing
to the limit $\gamma\to 0$, the macroscopic density profile
$\rho(r,t)$ satisfies the equation
\begin{equation}
\qquad \frac{\partial \rho}{\partial t}=\nabla   \cdot\left[\sigma_s(\rho)\nabla
\frac{\delta F(\rho)}{\delta \rho }\right]
\end{equation}
Here  $F(\rho)=\int f_s(\rho(r))\,dr- (1/2)\iint J(r-r')\rho(r)\rho(r')\,drdr'$,
where $f_s(\rho)$ is the (strictly convex) free energy density of the
reference system, and $\sigma_s(\rho)$ is the mobility of the system with only short-range
interactions.  See also \cite{GL2} and the review paper  \cite{GLP}.
The model is used to study phase segregation in \cite{GLM2000}.

Further motivations for  model \eqref{eq.ic}  can be found in \cite{CaVa09} and \cite{JLVAbel}, which contain references to applications in  dislocation dynamics and in superconductivity, as well as current mathematical progress.

\medskip

\noindent {\bf Mathematical results.} Paper \cite{CaVa09} contains the proof of existence of a weak solution of Problem \eqref{eq1}-\eqref{eq.ic} when  $u_0$ is a bounded function and has exponential decay at infinity. Besides, a  number of basic  properties are proved, like  energy estimates,  bounds in the $L^p$ spaces, and the property of finite propagation that says that compactly supported data produce solutions whose support is compact in space for every positive time. However, the  question of uniqueness of weak solutions is a pending open problem (in more than one space dimension). Comparison theorems, a crucial tool in parabolic equations, are only available under special circumstances (i.e., for so-called true super- or sub-solutions).  The asymptotic behaviour of the solutions as $t\to\infty$ has been studied by two of the authors in \cite{CaVa11} using obstacle problems and entropy estimates.

The regularity theory that we develop in the present paper is as follows:

a) If $u$ has initial data in $L^1\log L$, then it becomes instantaneously bounded;

b) If $u$ has initial data in $L^1$, then it falls into the previous case

c) Bounded solutions are continuous with a modulus of continuity.

After some preliminaries, Section \ref{sec.prelim}, and the needed theory on bilinear forms contained in Section \ref{sec.bf}, the  boundedness results are stated as Theorem \ref{th:L-inf} and proved in Section \ref{sec.reg}. The proof of the $C^\alpha$ regularity result, Theorem  \ref{mainthm}, takes up Sections \ref{sec.holder} to \ref{s-large}.  It says that for $0<s<1$, with $s\ne 1/2$, bounded solutions $u\ge 0 $ are $C^\alpha$ continuous in $(x,t)$ with some universal exponent $\alpha\in (0,1)$ that depends on $N,s$. The proof of this result in the range $0<s<1/2$ uses a number of techniques
that are becoming classical in the study of regularity of nonlocal diffusion problems, but it is complicated since we must take into account both the nonlinearity and the possible degeneracy. Section \ref{s-large}  covers the more difficult range $s>1/2$. The regularity result in that case uses transport ideas in the form of a geometrical transformation to absorb the uncontrolled growth of one of the integrals that appear in the iterated energy estimates. The case $s=1/2$ has new difficulties and will be treated separately.

As a consequence of these results, in Section \ref{sec.ext} we complete the existence theory by  constructing a continuous  weak solution  for any initial data $u_0\in L^1(\ren)$, $u_0\ge 0$.

\medskip

\noindent {\bf Notations.}
We will refer to Equation \eqref{eq1} as the FPME (with F for fractional). We will use the notation $(-\Delta)^{s}$ with $0<s<1$ for the fractional powers of the Laplace operator defined on the Schwartz class of functions in $\RR^N$ by Fourier transform and extended in a natural way to functions in the Sobolev space $H^{2s}(\RR^N)$.  Technical reasons imply that in one space dimension the restriction $s<1/2$ will be observed. The inverse operator is denoted by ${\cal L}_s=(-\Delta)^{-s}$ and can be realized by convolution
\begin{equation}
{\cal L}_s u=L_s\star u, \qquad L_s(x)=c(N,s)|x|^{-N+2s}.
\end{equation}
${\cal L}_s$ is a positive self-adjoint operator. We will write ${\cal H}_s={\cal L}_s^{1/2}$ which has kernel $L_{s/2}$. The subscript $s$ will be omitted when $s$ is fixed and known. For functions that depend on $x$ and $t$, convolution is applied for every fixed $t$ with respect to the space variables. We then use the abbreviated notation $u(t)=u(\cdot,t)$.

For a measurable  $u\ge 0$ and  for $k>0$ we denote by $u_k^+=(u-k)_+=\max\{u-k,0\}$, and $u_k^-=\min\{u-k,0\}$ in such a way that $u_k^+\ge 0\ge u_k^-$, the supports of $u_k^+$ and $u_k^-$ agree only on points where $u=0$, and also $u=k+u_k^+ + u_k^-$. We will use similar notations: $u_\varphi^+=(u-\varphi)^+$, $u_\varphi^-=(u-\varphi)^-$ when $\varphi$ is a function and not just a constant, and then we may split $u$ as follows: $u={\varphi}+u_{\varphi}^+ + u_{\varphi}^-$.

\section{Preliminaries. Existence and basic estimates}
\label{sec.prelim}

\noindent {\bf Definition.} {\sl We say that $u$ is a weak solution of Problem {\rm (\ref{eq1})-\eqref{eq.ic}} in \ $Q_T=\RR^N\times (0,T)$ \ with initial data \ $u_0\in L^1(\RR^N)$ if $u\in L^1(Q_T)$, ${\cal L}(u)\in L^1_{loc}(0,T:W^{1,1}_{loc}(\RR^N))$, and $u\,\nabla{\cal L}(u)\in L^1(Q_T)$,
and if the identity
\begin{equation}
\iint u\,(\eta_t-\nabla {\cal L}(u)\cdot\nabla\eta)\,dxdt+ \int
u_0(x)\,\eta(x,0)\,dx=0
\end{equation}
holds for all  continuously differentiable test functions $\eta$ in $Q_T$  that are compactly supported in the space variable and vanish near $t=T$.}

The following results have been proved in \cite{CaVa09}.

\begin{thm}\label{thm.ex} Let $u_0\in  L^\infty(\ren)$, $u_0\ge 0$, and such that
\begin{equation}
u_0(x)\le A\,e^{-a|x|} \qquad \mbox{for some $A,a>0$}\,.
\end{equation}
Then there exists a weak solution $u$ of Equation \eqref{eq1} with initial data $u_0$. Besides,  $u\in L^\infty(0,\infty: L^1(\ren))$, $u\in L^\infty(Q)$, $\nabla {\cal H}(u)\in L^2(Q)$. Moreover, for all $t>0$ we have conservation of mass:
\begin{equation}
\int_{\ren} u(x,t)\,dx=\int_{\ren} u_0(x)\,dx\,,
\end{equation}
as well as the $L^\infty$ bound: $\|u(t)\|_\infty \le \|u_0\|_\infty$. The constructed solution decays exponentially as $|x|\to\infty$. The first energy inequality holds in the form
\begin{equation}
\dint_0^t\dint_{\ren} |\nabla {\cal  H} u|^2\,dxdt +  \dint_{\ren}  u(t)\log(u(t))\,dx
 \le \dint_{\ren} u_0\log(u_0)\,dx\,,
\end{equation}
while the second says that for all $0<t_1<t_2<\infty$
\begin{equation}
\int_{t_1}^{t_2}\int_{\ren} u\,|\nabla {\cal L}u|^2\,dxdt+ \frac 12\int_{\ren} |{\cal H}u(t_2)|^2\,dx\le
\frac 12\int_{\ren} |{\cal H}(u(t_1)|^2\,dx\,.
\end{equation}
\end{thm}

\noindent {\bf Other properties of the constructed solutions.} Here are some of the most useful

\noindent $\bullet$ {Translation invariance.} The equation is invariant under translations in space and time, and this property reflects on the set of weak solutions.

\medskip

\noindent $\bullet$ {Scaling.} Moreover, the equation is invariant under a subgroup of the group of dilations in $(u,x,t)$, and this implies a scaling property for the set of solutions. Namely, if $u(x,t)$ is a weak solution as described in the existence theorem, with initial data $u_0(x)$, and $A,B,C$ are positive constants, then $\widehat u(x,t)= A\, u(Bx, Ct)$ is again a weak solution on the condition that $A=CB^{-2+2s}$. It has initial data $\widehat u_0(x)=Au_0(Bx)$.

\medskip

\noindent $\bullet$ Conservation of sign: $u_0\ge 0$ implies that $u(t)\ge 0$
for all times.

\medskip

\noindent $\bullet$ $L^p$ estimates. The $L^p$ norm of the solutions, $1< p\le \infty$, does not increase in time.

\medskip

\noindent $\bullet$ {Finite propagation:}
Compactly supported initial data $u_0(x)$ give rise to solutions $u(x,t)$ that have the same
property for all positive times, i.e., the support of $u(\cdot,t)$
is contained in a finite ball $B_{R(t)}(0)$ for any $t>0$.

 \medskip

\noindent  $\bullet$ A standard comparison result for parabolic equations does not work in general. This is one of the main technical difficulties in the study of this equation. In fact, special situations are found in \cite{CaVa09} where some comparison holds by using so-called true super- and sub-solutions.

\medskip

\noindent {\bf Energy solutions.} The constructed solutions are limits of smooth functions for which the energy inequalities are justified. In the sequel we will need this fact and also similar integrations by parts involved in the new energy inequalities.  In particular, we want the weak solution to satisfy the identities
 \begin{equation}
\iint u\,(\eta_t-{\cal B}_r(u,\eta))\,dxdt+ \int
u_0(x)\,\eta(x,0)\,dx=0,
\end{equation}
where ${\cal B}_r$ is the bilinear form that will be defined in the next section, $r=1-s$ and $\eta\in L^2(0,T:H^r(\ren))$, $\eta$  bounded, $\eta_t\in L^2(Q_T)$.  This class of solutions can be called weak energy solutions. Below (see Formulas \eqref{wes.f1}, \eqref{wes.f2}), we will need a version a this definition what consists in using $\eta= f(u)$ and integrating in time to get
\begin{equation}\label{wes1}
\left.\int  F(u(t))\,dx \right|_{t_1}^{t_2}+ \int_{t_1}^{t_2}\int \nabla \big[f(u)\big]\,u\nabla {\cal L}u\,dxdt=0.
\end{equation}
with  $f$ is smooth and bounded, $F(s)=\int^s f(s)\,ds$, and $0\le t_1<t_2\le T$. Note that we do not need to assume regularity for $\eta_t$. The constructed bounded solutions are energy solutions in this sense.
 We  will also use $\eta= f(u/\vp)$, where $f$ as before and $\vp(x)$ is a smooth positive function that does not vanish, see Formula  \eqref{wes.f3} and later.

\section{Bilinear forms}\label{sec.bf}

Before proceeding with the study of the boundedness and regularity properties, we need some results on fractional operators. The  bilinear form associated to the space $H^r(\ren)=W^{r,2}(\ren)$, $0<r<1$, is
\begin{equation}\label{eq.bilin}
{\cal B}_r(v,w)= C_N \diint (v(x)-v(y))\dfrac1{|x-y|^{N+2r}}(w(x)-w(y))\,dxdy
\end{equation}
It is easily seen by Fourier transform that this is well defined for functions in $H^r(\ren)$.
We will omit the subindex and write ${\cal B}$ instead of ${\cal B}_r$ when the context is clear.

\begin{cor} \label{cor3.1} (a) If $v$ is a monotone function of $w$, i.\,e., if $v=G(w)$ with $G'\ge 0$ then
$$
{\cal B}(v,w)\ge 0.
$$
\noindent (b) If $G'(s)\le C$ for some constant $C>0$, then
$$
{\cal B}(G(w),w)\le C\,{\cal B}(w,w).
$$
\end{cor}

\begin{prop}\label{prop.bilin} For every $u,v\in H^1(\ren)=W^{1,2}(\ren)$ we have
\begin{equation}\label{def.bilin}
{\cal B}_r(v,w)= \diint \nabla v(x)\dfrac1{|x-y|^{N-2+2r}}\nabla w(y)\,dxdy\,.
\end{equation}
\end{prop}

\noindent {\sl Proof.} Prove first for $C^\infty_c$ functions and use a smoothing and truncation of the kernel. Then pass to the limit. \qed

\noindent {\bf Remark.} Since the weak formulation of the FPME leads to an expression of this latter form with  kernel $K= c|x-y|^{-N+2s}$, we will put below $r=1-s$. Actually, all that we will use in  Sections \ref{sec.reg} and later, in accordance with Proposition \ref{prop.bilin}, is two kernels $L$ and $K$ such that $L, K\ge 0$ and $\Delta L =K$, as well as
the bilinear forms associated to the pairings
$$
\nabla (\cdot) \cdot L\cdot \nabla(\cdot) \longleftrightarrow \mbox{difference} \cdot K\cdot \mbox{difference},
$$
which is a short way of writing the equivalence of formulas \eqref{eq.bilin}, \eqref{def.bilin}.
In later calculations we will also use  the following positivity properties of the integration of these kernels applied to truncations of functions.

\begin{lem}\label{cor3.3} Let $u\in H^r(\ren)$ and $u_k^+\in H^r(\ren)$. Then,
\begin{equation}
{\cal B}_r(u_k^+,u)\ge {\cal B}_r(u_k^+,u_k^+).
\end{equation}
\end{lem}

\noindent {\sl Proof. } If  $K_r(x,y)=|x-y|^{-(N+2r)}$ is the kernel of ${\cal B}_r$, we have
$$
\begin{array}{l}
{\cal B}_r(u_k^+,u) - {\cal B}_r(u_k^+,u_k^+)={\cal B}_r(u_k^+, u-u_k^+)={\cal B}_r(u_k^+, k+ u_k^-)=\\
\diint (u_k^+(x)-u_k^+(y))\,K_r(x,y)\,(u_k^-(x)-u_k^-(y))\,dxdy\,.
\end{array}
$$
Now, given the fact that $u^+(x)u^-(x)=0$ a.e. and symmetry in $x,y$, the last integral equals
$$
-2\int u_k^+(x)u_k^-(y)\,K(x,y) \,dxdy\ge 0. \hskip 3cm \mbox{\qed}
$$

We will also need the following embedding inequality.

\begin{lem}\label{lem.q} For every $u\in L^1(\ren)\cap H^r(\ren)$ we have
\begin{equation} \label{lem3.4}
\int u^q\,dx \le C\|u\|_{L^1}^\theta\|u\|_{H^r}^2\,,
\end{equation}
where  $\theta=2r/N$, $q=2+\theta=2+(2r/N)$, and $C=C(N,r)>0$. Moreover, for  every $u\in L^2(\ren)\cap H^r(\ren)$ we have
\begin{equation} \label{lem3.4b}
\int u^{q_1}\,dx \le C\|u\|_{L^2}^\theta\|u\|_{H^r}^2\,,
\end{equation}
 where $\theta=4r/N$, $q_1=2+(4r/N)$, and $C=C(N,r)>0$  as before.
\end{lem}

\noindent {\sl Sketch of the proof.} We use the Sobolev inequality that says that
\begin{equation} \label{sob}
\left(\int u^p\,dx\right)^{2/p}\le C \|u\|^2_{H^{r}}
\end{equation}
for some $p>2$  depending on $r\in (0,1)$ and $N$. $C$ depends also on $r$ and $N$.
Actually, $p=2N/(N-2r)$, when $N\ge 2$ or when $N=1$ if, in addition, $ 0<r<1/2$. We want to control $\|u(t)\|_q$ for some $q>1$ using that
$$
\int u^q\,dx \le \left(\int u\,dx\right)^{\theta}\left(\int u^p \,dx\right)^{1-\theta},
$$
where $q=\theta+p(1-\theta)$. We will take the values
$$
\theta=\frac{p-2}p, \quad 1-\theta=\frac2{p}.
$$
The proof when $u\in L^2$ is quite similar.

The case $N=1$, with $1/2\le r<1$ is easy. Take $0<r'<1/2$ and observe that $L^1\cap H^{r}\subset  H^{r'}$ as a continuous embedding.
Now we can use (\ref{sob}), with $r'$ replacing $r$, and get (\ref{lem3.4}) with, perhaps, some different values of $q$ and $\theta$. \qed

\medskip

\begin{rem} \label{convex}
We also recall that for every convex function $\Phi$ the quantity $\int \Phi(u(t))\,dx$ is non-increasing in time in the FPME evolution. We note that
for such a function $\Phi$ we have
$$
\left.\int \Phi(u(x,t))\,dx\right|_{T_0}^{T_1} + \int_{T_0}^{T_1}{\cal B}(G(u),u)\,dt \le 0\,,
$$
where $G(u)$ is the primitive of $\Phi''(u)u$. The bilinear form ${\cal B}$ is as above and
$$
{\cal B}(G(u),u)\ge c\,{\cal B}(u,u)\sim \|u\\|^2_{H^r}
$$
if $G'$ is strictly positive. We will use this in the case of the truncations  in the form $\Phi(u_k^+)$ with $u_k^+=(u-k)^+$ and then  $G=0$ for $u \le k$.
\end{rem}

%
\section{Boundedness of solutions}\label{sec.reg}

 This section is devoted to proving the main boundedness result. Here $0<s<1$.

\begin{thm}\label{th:L-inf} Let $u$ be a weak energy solution of Problem \eqref{eq1}--\eqref{eq.ic} with $u_0\in L^1(\mathbb{R}^N)
\cap L^\infty(\mathbb{R}^N)$,  as constructed in Section 2. Then, there exists a positive constant $C$ such that for
every $t>0$
\begin{equation}\label{form:L-inf}
\sup_{x\in\mathbb{R}^n}|u(x,t)|\le C\,t^{-\alpha }\|u_0\|_{L^1(\mathbb{R}^n)}^{\gamma}
\end{equation}
with precise exponents $\alpha=N/(N+2-2s)$, $\gamma=(2-2s)/(N+2-2s)$. The constant $C$ depends only on $N\ge 1$ and $s\in (0,1)$.
\end{thm}

In dimension $N=1$ this is proved in \cite{BKM}. Our proof  applies to all $N\ge 1$ and is divided into three subsections. Note that this estimate and conservation of mass imply a decay for all intermediate norms $L^p$ with $1<p<\infty$:
\begin{equation}\label{form:L-p}
\|u(\cdot,t)\|_p\le C_p\,t^{-\alpha_p}\|u_0\|_{L^1(\mathbb{R}^n)}^{\gamma_p}\,,
\end{equation}
where $\alpha_p=\alpha(p-1)/p$ \ and $\gamma_p=(1+\gamma(p-1))/p$.

\subsection{Better  integrability properties for solutions}

We start with the following partial result.

\begin{lem} {\rm  (\lq\lq From $L^1$ to $L \log L$")}  Let $u\ge 0$ be a weak energy solution  of the FPME. If the initial data are integrable, then $u(\cdot,t)\in L \log L$ for all positive $t>0$ and for all small $t\le t_0$ we have
\begin{equation}
\int u(t)\log(1+u(t))\,dx +  \frac{1}{t}\int_0^t s \,{\cal B}(u(s)_1^+,u(s)_1^+)\,ds \le C_0  |\log(t)|\|u_0\|_1\,,
\end{equation}
where $C_0$ depends only on $N,s$, and we write $u(t)_1^+=(u(t)-1)_+$. The time $t_0$ is estimated as $t_0=\inf\{1, C\|u_0\|_{L^1}^{-\vartheta}\}$ for some constants $C,\vartheta>0$. \end{lem}

\noindent {\sl Proof.}   We use as test function $\eta=t\log(1+ u)$ on the weak form of the equation.
 After observing that  $h(u)=(u+1)\log(1+u)\ge 0$  satisfies $\log(1+u)=h'(u)-1$ and  integrating in $\ren\times (0,\tau)$, we get the identity
\begin{equation}\label{wes.f1}
\tau\hskip-4pt\int h(u(\tau))\,dx+\int_0^\tau t {\cal B}(g(u),u)\,dt=\int_0^\tau\int h(u(t))dxdt+\tau\hskip-4pt\int u(\tau)\,dx-\iint u\,dxdt\,,
\end{equation}
  where ${\cal B}={\cal B}_r$ with $r=1-s$, as already explained, and $g'(u)=u/(1+u)$ with $g(0)=0$.   The last two terms in the display disappear by mass conservation. Note next that  $g'(u)\ge 1/2$   for $u\ge 1$, hence $g(u)\ge (1/2)(u-1)^+$, so that, writing $u_1^+=(u-1)^+$ we have by the already
  monotonicity properties of $\cal B$:
  $$
  {\cal B}(g(u),u)\ge (1/2){\cal B}(u_1^+,u)
  \ge (1/2){\cal B}(u_1^+,u_1^+).
  $$
With this we arrive at
\begin{equation}\label{ineq.te}
\tau\int h(u(\tau))\,dx+\frac12 \int_0^\tau t {\cal B}(u_1^+,u_1^+)\,dt\le \int_0^\tau\int h(u(t))dxdt\,.
\end{equation}
Recalling the definition of $h$ and using $\log(1+u)\le u$, we get
$$
\sup_{\{0<t<\tau \}} t\int u(t)\log(1+u(t))\,dx + \frac12 \int_0^\tau t {\cal B}(u_1^+,u_1^+)\,dt\le \int_0^\tau\!\!\!\int u\log(u+1)\,dxdt+ \tau \int u_0\,dx.
$$

We still have  a ``bad term'' in the right-hand side containing $\int u\log (1+u)\,dx$, and it  is tackled as follows: note that the left-hand side controls \  $\int tdt \int u^q\,dx$ \ if  $q$ is as in Lemma \ref{lem.q}; take such a $q>2$. Then, for any $M>2$ we have
$$
\int u\,\log(1+u)\,dx \le \log (M+1)\int_{u<M} u\,dx + \frac{\log (1+M)}{M^{q-1}}\int_{u>M} u^q\,dx\,,
$$
where we have used the fact that $u\log(1+u)/u^q$ is decreasing in $u$ for $u\ge M$. Next,  there is a constant $C(N,q)>0$ such that $u^q\le C(u-1)_+^q$ for this range of $u$.
Choose now $M=t^{-\alpha}$ and then $\beta<\alpha (q-1)-\ve$ (taking care that $\beta>1$). With all this, we get for small $\tau$
$$
\begin{array}{c}
\dint_0^\tau\dint u\,\log(1+u)\,dxdt \le C_1 \dint_0^\tau dt |\log(t)|\big(\sup_t \int u \big)+ C_2\dint_0^\tau  t^{\beta} dt (\int (u-1)_+^{q}\,dx )
 \\[10pt]
\le C_1 \tau|\log(\tau)|\dint u_0\,dx + C_2\dint_0^\tau t^{\beta} \|u(t)_1^+\|_{L^1}^\theta\|u_1^+(t)\|_{H^r}^2\,dt\,.
\end{array}
$$
Since $\beta>1$ the last term is controlled by the ${\cal B}$-energy term in the left-hand side of \eqref{ineq.te} for small $\tau$. In particular, we choose $0<\tau<1$ and $\tau^{\beta-1}\le (2C_2\|u_0\|_{L^1}^{\theta})^{-1}$. The other term is a multiple of $\int u_0$, hence bounded. We get
$$
\sup_{\{0<t<\tau \}} t\int u(t)\log(1+u(t))\,dx +  \int_0^\tau t {\cal B}(u_1^+,u_1^+)\,dt\le C_3\tau  |\log(\tau)|\,\|u_0\|_1.
$$
Putting $t=\tau$ we get the result.  \qed

\begin{lem}  {\rm  (\lq\lq From $L \log L$ to $L^2$")}   Initial data in the space $L \log L$ imply that  $u(\cdot,t)\in L^2$ for all positive $t>0$ and the bound on the $L^2$ norm of $(u-1)^+$ depends only on $t$, $s,N$,  $\|u_0\|_1$, and $\|(u_0-1)^+\|_{L\log L}$.
\end{lem}

\noindent {\sl Proof.}  We define $v(x,t)=u(x,t) \vee 1$ so that $v=1+u_1^+ \ge 1$ and $v_t=u_t\,\chi(u>1)$. Recall the notation $u_k^+=(u-k)^+$. Then,
$$
\frac{d}{dt}\int [v\log(v)-v]\,dx=\int \log(v) \, v_t= \int_{u>1} \log(v) \,u_t.
$$
Using the weak form of the equation with $\eta=\log v$ as  test function,   we  get
\begin{equation}\label{wes.f2}
\left.\int  \,[v\log(v)-v](t)\,dx \right|_{t_1}^{t_2}+ \int_{t_1}^{t_2}\int \nabla \big[\log(v)\big]\,u\nabla {\cal L}u\,dxdt=0.
\end{equation}
We work out the last term for fixed time and observe that, since $u=v$ for $u>1$, we have
$$
\int_{u>1} \nabla \big[\log(v)\big]\,u\nabla {\cal L}u\,dx=\int \nabla v\cdot \nabla {\cal L}u\,dx=
\int \nabla u_1^+ \cdot \nabla {\cal L}u\,dx={\cal B}_r(u_1^+,u)\,.
$$
Using again the monotonicity of $\cal B$, see Corollary \ref{cor3.1} and Remark \ref{convex}, and putting  $h(v)=v\log(v)-v$,  we get
\begin{equation}
\int  \,h(v(t_2))\,dx + \int_{t_1}^{t_2} {\cal B}_r(u_1^+,u_1^+)\,dt \le
\int  \,h(v(t_1))\,dx\,.
\end{equation}
Note that  $h(v) $ is convex for $v>1$  and the right-hand side  is bounded by a combination of $\|u_0\|_1$ and $\|(u_0-1)^+\|_{L\log L}$. Hence,
\begin{equation} \label{quan}
 \int_{0}^{t_2} {\cal B}_r(u_1^+,u_1^+)\,dt \le C.
 \end{equation}
 Recall finally that $ {\cal B}_r(u_1^+,u_1^+)\sim \|u_1^+(\cdot,t)\|^2_{H^r}, $ with $r=1-s$. Use Lemma \ref{lem.q} to conclude that
$u_1^+(t)\in L^q(\ren)$ for some $q>2$. More quantitatively, this together with (\ref{quan}) and Remark \ref{convex} give the estimate
\begin{equation}
\sup_{t>0} t\|u_1^+(t)\|^q_{L^q}\le C\,,
\end{equation}
with $C$ depending as in the statement of the lemma. Interpolation with $L^1$ gives the result. \qed

\subsection{Boundedness}\label{subs.endbound}

With the preceding results, we may assume that  $u_0\in L^1(\ren)\cap L^2(\ren)$ after some displacement of the time origin. Then we can follow the De Giorgi approach (as outlined for instance in \cite{CVass-DG},  \cite{CVass}). We consider the   truncations
$$
  u_j^+(x,t)=(u-M(2-2^{-j}))_+
$$
 The value of constant $M>0$ will be conveniently chosen later. Actually, we may assume $\int u_0^2\,dx$ very small by selecting $M$ large.

 \noindent {\bf Claim.} The following energy inequality holds  for all $0<t_1<t_2$:
\begin{equation} \label{bdd}
\int    u_j^+(t_2)^2\,dx + 2M(2-2^{-j})\int_{t_1}^{t_2} {\cal B}(  u_j^+,  u_j^+)\,dt \le
\int    u_j^+(t_1)^2\,dx\,.
\end{equation}
To see this, we use the definition of weak solution for our FPME with $\eta=  u_j^+$ as a test function. Then, for $t_1<t_2$, we have
$$
\int_{t_1}^{t_2}\int  u_t   u_j^+ \,dx\, dt= \frac 12 \int    u_j^+(t_2)^2\,dx - \frac 12 \int    u_j^+(t_2)^2\,dx.
$$
For the RHS we observe that $u=  u_j^++M(2-2^{-j})$ whenever $u\ge M(2-2^{-j})$. Hence,
$$
\int_{u>M(2-2^{-j})}u\nabla   u_j^+ \cdot \nabla {\cal L}u\,dx=\frac 12 \int \nabla (  u_j^+)^2 \cdot \nabla {\cal L}u\,dx+ M(2-2^{-j}) \int \nabla   u_j^+ \cdot \nabla {\cal L}u\,dx =
$$
$$
=\frac 12 {\cal B}((  u_j^+)^2,u) + M(2-2^{-j})\,{\cal B}(  u_j^+,u)\ge  M(2-2^{-j})\,{\cal B}(  u_j^+,  u_j^+),
$$
In the last inequality we have used both Corollary \ref{cor3.1} and Lemma \ref{cor3.3}. This gives (\ref{bdd}). \qed

\medskip

We now fix $t_0>0$. We want to prove that the solution $u$ is bounded for all times $t\ge t_0$. As in \cite{CVass}, let us define the total energy for the truncated function $u_j^+$ as
\begin{equation*}
{\cal A}_j=\sup_{t\ge T_j}\int  (u_j^+)^2(t)\,dx + 2M\int_{T_j}^\infty {\cal B}( u_j^+, u_j^+)\,dt,
\end{equation*}
where $T_j=t_0(1-2^{-j})$. From (\ref{bdd}), taking arbitrary values $t_2=t\ge T_j$ and $t_1=t' \in [T_{j-1},T_j]$ we have
\begin{equation} \label{bdd2}
{\cal A}_j\leq \inf_{t' \in [T_{j-1},T_j]} \int  (u_{j}^+)^2(t')\,dx.
\end{equation}
Observe now that $u_{j}^+(x)>0$ implies $u_{j-1}^+(x)=u_{j}^+(x)+M2^{-j}>M2^{-j}$. Therefore, for every $p>2$ we have (keeping the time fixed)
$$
\int (u_{j}^+)^2\dx =\int (u_{j-1}^+ -M2^{-j})^2\cdot \chi_{\{u_{j}^+>0\}}\, dx
$$
$$
 \leq
\int (u_{j-1}^+)^2\left(\frac {u_{j-1}^+}{2^{-j}M}\right)^{p-2}\, dx=\left(\frac {2^j}{M}\right)^{p-2}\dint (u_{j-1}^+)^p\, dx.
$$
If $p>2$ is the exponent corresponding to Sobolev's embedding theorem, we deduce from (\ref{bdd2}) that
$$
\begin{array}{l}
{\cal A}_j\leq \inf_{t' \in [T_{j-1},T_j]} \left(\dfrac {2^j}{M}\right)^{p-2}\dint (u_{j-1}^+(t'))^p\, dx\\[10pt]
\le C_N \left({2^j}/{M}\right)^{p-2}  \inf_{t' \in [T_{j-1},T_j]}\left[{\cal B}( u_{j-1}^+, u_{j-1}^+)(t')\right]^{p/2}.
\end{array}
$$
Taking averages in $t'$ we arrive to the inequality
$$
{\cal A}j\leq C_N
\left({2^j}/{M}\right)^{p-2}  \left[\frac 1{T_j-T_{j-1}}\int_{T_{j-1}}^{T_j}{\cal B}( u_{j-1}^+, u_{j-1}^+)(t')\,dt' \right]^{p/2}
$$
$$ \leq C_N
\frac {2^{j(p-2)}2^{jp/2}}{M^{p-2} (Mt_0)^{p/2}} \left[M\int_{T_{j-1}}^{\infty}{\cal B}( u_{j-1}^+, u_{j-1}^+)(t')\,dt'\right]^{p/2}.
$$
This leads to a recurrence relation  of the form
$$
{\cal A}_{j}\le C^j \left(\frac {{\cal A}_{j-1}}{Mt_0}\right)^{1+\delta},
$$
with $\delta=\frac p2 -1>0$, that implies ${\cal A}_{\infty}=0$ if ${\cal A}_0/Mt_0$ is small. This determines the correct value of $M$ to choose. The conclusion is then that $u(x,t)\le 2M$ for all $t\ge t_0$. \qed

\subsection{End of proof of the theorem}

The preceding subsections have established the result for any fixed $t>0$, and we know that
\begin{equation}
u(x,t)\le C(N,s,\|u_0\|_1,t)
\end{equation}
but we do not know the dependence of $C$ on its arguments in a precise way. We need to prove that this dependence takes the form \eqref{form:L-inf}. This is just a consequence of the scaling group that allows to pass from a solution $u(x,t)$ to the rescaled solution
\begin{equation}
{\widetilde u}(x,t)=A\, u(L\, x,T\,t)
\end{equation}
on the condition that $A=TL^{-2+2s}$. On the other hand, we want to reduce  $\widetilde u$ to unit mass, $\int \widetilde u(x,t)\,dx=1$, and this means $A=L^N/\|u_0\|_1$. All together this gives (with $\|u_0\|_1=M$)
$$
L=(MT)^{\beta} \quad A=M^{2(1-s)\beta}T^{-N\beta},
$$
where $\beta=1/(N+2-2s)$. We now apply the boundedness result to $\widetilde u$ at $t=1$,
i.\,e.,
$$
\widetilde u(x,1)\le C(N,s) \qquad \forall x\in \ren.
$$
Going back to the $u$, we have
$$
\sup_x u(x,T)= A\,\sup_x {\widetilde u}(x,1)=M^{2(1-s)\beta}T^{-N\beta}\,C(N,s)
$$
which gives the desired result upon replacing $T$ by $t$. For other instances of this scaling argument cf. \cite{JLVSmoothing}. \qed


\section{H\"older regularity. Main result and basic lemmas}\label{sec.holder}

Once the question of boundedness is settled, we proceed with the local regularity of solutions.
This is the main result.

\begin{thm}\label{mainthm} Let $u\ge 0$ be a bounded weak energy solution defined in a space-time strip $S=\RR^N\times [T_1,T_2]\subset \re^{N+1}$. Let $s\in (0,1)$, $s\ne 1/2$. Then  $u$ is $C^\alpha$ continuous in the interior of $Q$
for some exponent $\alpha(N,s)\in (0,1)$  and a constant that depends also on the dimensions of the subdomain $Q\Subset S$ and the bounds on $u$ in $L^\infty(S)$ and $L^\infty((T_1,T_2):L^1(\RR^N))$.
\end{thm}

\noindent {\bf Strategy.}  Since the equation is space- and time-invariant we may assume that $T_1<T_2 =0$, and then we may study the regularity around $x=0$ and $t=0$.

The main ideas are two: on the one hand, we will prove some basic De Giorgi-type oscillation lemmas that say that the oscillation of the solution $u$ decays when we restrict a basic domain, say, the cylinder $\Gamma_4=[-4,0]\times B_4(0)$, into a smaller cylinder like $\Gamma_1=[-1,0]\times B_1(0)$. The second ingredient is the scaling property of the equation that allows to renormalize the solution through the  transformation
\begin{equation}\label{scaling}
\widehat u(x,t)= A\, u(Bx, Ct)
\end{equation}
with $C,B>0$ free parameters, and  $A=CB^{-2+2s}$. The way of attacking the problem is through the iterated application of the lemmas. At the end of every step we renormalize the solution defined in $\Gamma_1$ into a rescaled solution defined in $\Gamma_4$ and we start a new application of the oscillation lemmas. In this way, we will show that the oscillation of the solution $u$ decays dyadically in a family of space-time cylinders shrinking dyadically to a point.

The needed lemmas have a simpler expression for $0<s<1/2$ where the diffusion is more similar to the
standard porous medium case. For $1/2\le s < 1$  convection effects appear that make some integrals diverge, and this makes the analysis more difficult, needing new techniques. The detailed study of how to proceed in the case $s>1/2$ are contained in Section \ref{s-large}. Until then we assume that $s<1/2$

\medskip

\subsection{The oscillation reduction lemmas}

These technical results need only be proved for bounded nonnegative weak solutions defined in a strip $S_R=[-R,0]\times \RR^N$. We denote by $\Gamma_R$ the parabolic cylinder $[-R,0]\times B_R(0)$. One of the lemmas controls the decrease  of the supremum  of the solution once we restrict the size of the parabolic neighborhood of $(0,0)$, the other one implies that under suitable assumptions the solution separates from zero. A third one improves the first result so as to obtain a real alternative between going a bit down and a bit up. which leads to the proof of regularity. Here is the first  basic lemma.

\begin{lem} \label{reg.1} Let $0<s<1/2$. Given $\mu\in (0,1/2)$ and  $\ve_0$ small enough $(\mbox{in particular,} \ve_0\le 1-2s)$, there exists $\delta>0$ (depending possibly on $\mu,\ve_0, s$, and $N$) such that if we assume that \\ {\rm (i)}  the solution $u$ is bounded above in the strip $S_4=\re^n\times [-4,0]$ by
\begin{equation}\label{est.apriori1}
 \overline{\Psi}(x)= 1+(|x|^\ve-2)_+, \qquad 0<\ve <\ve_0\,,
\end{equation}
 and {\rm (ii)} $u$ is mostly below the level  $1/2$  in $\Gamma_4=B_4(0)\times [-4,0]$ in the sense that
\begin{equation}
|\{u>1/2\}\cap \Gamma_4|\le \delta |\Gamma_4|\,,
\end{equation}
 then we can lower the upper bound inside a smaller cylinder in the following quantitative form:
$ \left. u\right|_{\Gamma_1}\le 1-\mu$.\end{lem}

We summarize the result by saying that {\sl ``being mostly below 1/2 in space-time measure pulls down the supremum in a smaller nested cylinder''}. Note that for this lemma $\delta$ can be chosen as a non-increasing function of $\mu$ with $\delta(1/2+)=0$.  Also, $\ve$ can be chosen as small as we want by sacrificing the gain in oscillation. We also remark that the size of the cylinders can be changed, though this affects the values of $\delta$ if the new sizes do not conform with the parabolic scaling. Finally, the levels $u=1/2$ and $u=1$ are taken
by convenience, any pair of levels $0<M_1<M_2$ will do, though in principle the value of $\delta$ will change.

A similar result applies from below  but the proof is different since the equation is degenerate at $u=0$. The idea is that if $u$ is very often far from zero in $\Gamma_4$ then in a smaller, suitably nested cylinder  $u$ stays uniformly away from zero. The technical version explains how ``most of the space-time above $1/2$, pulls up away from zero''.

\begin{lem} \label{reg.1b}  Under the  same assumptions, given $\mu_0\in (0,1/2)$ there exists $\delta>0$ (depending possibly on $\mu_0,\ve_0, s$, and $N$) such that if  $u$ satisfies
\begin{equation}
|\{u \ge 1/2\}\cap \Gamma_4|\ge (1- \delta) |\Gamma_4| \,,
\end{equation}
 then $ \left. u\right|_{\Gamma_1}\ge \mu_0$.
\end{lem}

Again, $\delta$ is a non-increasing function of $\mu_0$. Let $\delta_0=\delta(1/4)$, that is, when $\mu_0=1/4$.
A more elaborate version of this lemma will be needed in some cases of the iteration.

We complement these two lemmas with a lemma that replaces the sentence \lq\lq most of the space-time" of Lemma \ref{reg.1} by \lq\lq in some set of positive measure".

\begin{lem} \label{reg.2} {\rm (\lq\lq Some of the space-time below, pulls down'')} Assume as before that $0<s<1/2$ and  $u$ is trapped between $0$ and  $\overline{\Psi}$ in $S_4$. Besides, assume now that
\begin{equation}
|\{u<1/2\}\cap \Gamma_4|\ge \delta_0  |\Gamma_4|,
\end{equation}
with $\delta_0$ defined as above. Then $ \left. u\right|_{\Gamma_1}\le 1-\mu'$,  for some  $\mu'(\delta_0)$.
\end{lem}

Notice that this second lemma applies only in one direction, reducing the oscillation from above.
As in the classical porous media, we cannot expect this lemma to hold in the \lq\lq pulling-up" case, due to the property of finite propagation (existence of solutions with compact support), a consequence of the degeneration of the equation. Nevertheless, this one-sided improvement will be enough  to prove that the oscillation decays dyadically as follows:

The iterated  use of the Lemma \ref{reg.2} from above after rescaling at every step, as long as it is possible, reduces the oscillation of $u$ from above and we start iterating and renormalizing to get $u_2, u_3,\cdots$. We note though that, as we do that, a renormalized solutions develop a ``tail" in the sense that the functions $u_k$ start to grow at infinity by effect of the scaling (both in vertical and horizontal directions). This is the reason for the form of the upper barrier that we use, which has an  $\ve$-tail.  Indeed, after $k$ steps, $u$ will be bounded by $(1-\mu)^{-k}$  outside the $B^{-k}$ dilation using \eqref{scaling}.  This first difficulty can be dealt by playing the integrability of the kernels $L$ or $K$ at infinity against a slow power growth  in $u$. Indeed, by sacrificing the modulus of H\"older continuity we may assume that the gain is very tiny  (i.e., $\mu$ very small). Then, the build up in $u(y,t)$ as $y$ tends to infinity will be very slow (like $|y|^\ve$), being absorbed by $K$. If Lemma \ref{reg.2} never fails along the iteration, then we are at a point where $u=0$ and a H\"older exponent is also found.

On the other hand, if the process breaks down, then the first time that Lemma \ref{reg.2} fails, it puts us into the hypotheses of Lemma \ref{reg.1b} and that pulls $u$ away from zero by a fixed amount $\mu_0$ (the need for an alternative is what makes Lemma  \ref{reg.1} insufficient). Then the operator becomes non-degenerate in the subsequent iterations and a counterpart of Lemma \ref{reg.2} applies also from below (upwards) since $u$ will always be bounded between $\mu_0$ and $1$.
From there on the gain on the oscillation of $u$ may come from above or below in the dyadic way we have shown. The details of the iteration are given in Section \ref{sec.finalCalpha} after the lengthy and delicate proof of the Lemmas is completed. We recall that all this will be done below for $0<s<1/2$.

\section{Lowering the maximum. Proof of  Lemma \ref{reg.1} }\label{sec51}

We start here the technical work. The basic idea in  the proof of the result is a particular kind of  ``localized energy inequalities'' that will be iterated
in the style of Giorgi to obtain the reduction on the maximum in a smaller domain. Localization is obtained by using a suitable sequence of cutoff functions.  In order to deduce the necessary energy inequalities we use integration by parts formulas and analysis of the  kernels.
A main role is played by the bilinear form ${\cal B}(v,w)$ as defined in \eqref{prop.bilin} with kernel $K(x)=c|x|^{-(N+2r)} $ and $r=1-s$. Moreover, we put $L(x)=c_1|x|^{-N+2s}$ so that $-\Delta L =K$. We will repeatedly use the equivalent form justified by Proposition \ref{prop.bilin}.

\begin{lem} Let $u$ be smooth and  compactly supported, let $v$ have small growth at infinity. Then,
\begin{equation}
{\cal B}(u,v)=  \diint (u(x)-u(y))K (x-y) (v(y)-v(x))\,dxdy.
\end{equation}
\end{lem}

We take a weak energy solution defined in a strip $Q=(-T,0) \times \ren$ in the sense of Section \ref{sec.prelim}. We justify the computations by recalling that $u$ can be approximated by smooth positive solutions of similar problems as done
in \cite{CaVa09}.

\medskip

\noindent {\bf \ref{sec51}.1.   An energy formula}.  We  consider a sequence of cutoffs $\vp(x)$ that have the form of perturbations of the level $u=1$  within a region containing the unit ball $B_1(0)$,  and an  ``outer wing'' rising up above the 1-level for larger values of $|x|$.  An explicit choice suiting our purposes will be done below. We only need to know at this stage that the cutoff function $\vp$  is smooth, lies above $ 1/2$ everywhere, and also that $u\le \vp$ for all $|x|\ge 3$ for all times $-4\le t\le 0$.

We use the function $\eta=\log((u/\vp)\vee 1)=\log(g)$ as a  test function in the weak form of the equation. Note that
$$
g:= \frac{u}{\vp}\vee 1 = 1+\frac{(u-\vp)_+}{\vp}=1+\frac{u_\vp^+}{\vp},
$$
where $u_\vp^+=(u-\vp)^+$ according to the adopted notation. Note that $g\ge 1$ and $g>1$ iff $u>\vp$. According to our assumptions,  $u_\vp^+$ and $g-1$ have compact support in the ball of radius 3.  We will often split $u$ as follows
$$
u=u_\vp^++\vp + (u-\vp)^-
$$
where we write $(u-\vp)^-=(u-\vp)\wedge 0= u_\vp^-$. Notice that with this notation we have  $u_\vp^- \le 0.$ After applying  the weak formulation of the equation with $\eta$ as above, we get on the LHS
for $T_1\le t\le T_2\le 0$:
\begin{equation}\label{wes.f3}
\left\{\begin{array}{l}
\dint_{T_1}^{T_2}\dint \eta\,\partial_t u\,dxdt =
\dint_{T_1}^{T_2}\dint \log\left(\frac{u}{\vp}\vee 1 \right)\vp\,\partial_t(u/\vp)\,dxdt=\\[18pt]
\dint \vp \left(\frac{u}{\vp}\vee 1 \right)
\left[ \log\left(\frac{u}{\vp}\vee 1\right)-1\right] \,dx \Big|_{T_1}^{T_2}=
\dint \vp (g\,\log g -g) \,dx \Big|_{T_1}^{T_2}\,.
\end{array}\right.
\end{equation}
We  will need an estimate of this quantity: after  adding 1 to the last integrand we get  the expression $H(g):=g\,\log g +1-g$ for which we have  the estimate $H(g)\sim (g-1)^2$
 for for $ 1\le g \le 2$, in the sense that
\begin{equation}\label{ineq.H1}
\frac12\,\left(\frac {u_\vp^+}{\vp}\right)^2 \le H(g)\le \left(\frac {u_\vp^+}{\vp}\right)^2\,.
\end{equation}

Let us now calculate  the right-hand side of the expression in the  weak formulation of the equation.
We have
\begin{equation}\label{star}
\begin{array}{l}
\dint dt \dint \log(g(x))\,\, \mbox{\rm div} \left[ u(x)\, \nabla L (x-y)\, (u(y)-u(x))\right] \,dxdy\\
=-\dint dt \dint_{u>\vp}\frac{\nabla g(x)}{g(x)} \, u(x) \nabla_x L (x-y)\, (u(y)-u(x)) \,dxdy= {\bf I}+{\bf II}\,,
\end{array}
\end{equation}
where we pass from the first line to the second integrating by parts. Recalling that
$$
g= \frac {u}{\vp}\vee 1 = 1+ \frac{u_\vp^+}{\vp}, \qquad \nabla g=\nabla(u_\vp^+/\vp)\,\chi(\{u>\vp\})\,,
$$
the first part of the  splitting is:
$$
\begin{array}{l}
{\bf I}=-\dint dt \diint u_\vp^+(x)(-\Delta L)(x-y)[u(y)-u(x)]\,dxdy\\[12 pt]
=-\dint dt \diint u_\vp^+(x)\,K(x-y)[u(x)-u(y)]\,dxdy\,.
\end{array}
$$
After symmetrizing, we get ${\bf I}=- \frac12 \int  {\cal B}(u_\vp^+,u)\,dt$, where
$$
{\cal B}(u_\vp^+,u):=\diint (u_\vp^+(x)-u_\vp^+(y))\, K(x-y) (u(x)-u(y))\,dxdy.
$$
On the other hand,
$$
{\bf II}=\dint dt\diint u_\vp^+(x)\frac{\nabla \vp(x)}{\vp(x)} \nabla_x L(x-y)[u(y)-u(x)]dxdy : =\int  {\cal Q}(u_\vp^+,u)\,dt.
$$

\noindent $\bullet$ In order to separate the good and bad components of both ${\cal B}$ and ${\cal Q}$, we use the decomposition $u=u_\vp^++\vp+ u_\vp^-$. We get
$$
{\cal B}(u_\vp^+,u)= {\cal B}(u_\vp^+,u_\vp^+)+ {\cal B}(u_\vp^+,\vp)+ {\cal B}(u_\vp^+, u_\vp^-).
$$
and a similar expression for $\cal Q$. We now make some observations:

 (i) ${\cal B}(u_\vp^+,u_\vp^+)$ is a positive quadratic form. We will pass the corresponding part of $\bf I$ to the LHS as a term with positive sign and thus complete the energy expression in the energy inequality that we want to derive.

(ii) ${\cal B}(u_\vp^+,u_\vp^-)$ has also the correct nonnegative sign because of these facts: $u_\vp^+$ and $u_\vp^-$ have opposite signs and disjoint supports, and $K\ge 0$. We could drop this term in a first calculation, but we will keep it and use it to control some of the bad terms in  $\cal Q$.

Summing up, we have up to now the basic identity for $T_1<T_2\le 0$:
\begin{equation}\label{energ.id}
\left\{\begin{array}{l}
\dint \vp (g\,\log g+1 -g)  \,|{_{T_2}} \,dx + \frac12\dint_{T_1}^{T_2}{\cal B}(u_\vp^+,u_\vp^+)\,dt+
\dfrac12 \dint_{T_1}^{T_2}{\cal B}(u_\vp^+,u_\vp^-)\,dt\\[10pt]
= \dint \vp (g\,\log g+1 -g) \,|{_{T_1}} \,dx -\dfrac12 \dint_{T_1}^{T_2}{\cal B}(u_\vp^+,\varphi)\,dt
+ \dint_{T_1}^{T_2}  {\cal Q}(u_\vp^+,u)\,dt\,.
\end{array}\right.
\end{equation}

(iii) We will think of the LHS as the basic energy of this calculation, and the RHS as the terms still to be controlled.

\medskip


\noindent {\bf \ref{sec51}.3. Cutoff functions, control of the RHS and final goal}

In order to tackle the RHS and continue the proof of the lemma, we need to make a convenient choice of the sequence of cutoffs. Though only some simple bounds on the functions and their derivatives  are used, a possible practical choice is as follows:
 \begin{equation}
\vp_k(x)=\min\{1+(|x|^\ve-2)_+, \ \overline{\vp}_k(x)\}, \qquad \overline{\vp}_k(x) = \frac78 + \frac{|x|^2}{16}-\frac12 4^{-k}\,,
\end{equation}
for some small $\ve>0$ and $k=1,2,\dots$

\begin{center}
\includegraphics[width=8.5 cm]{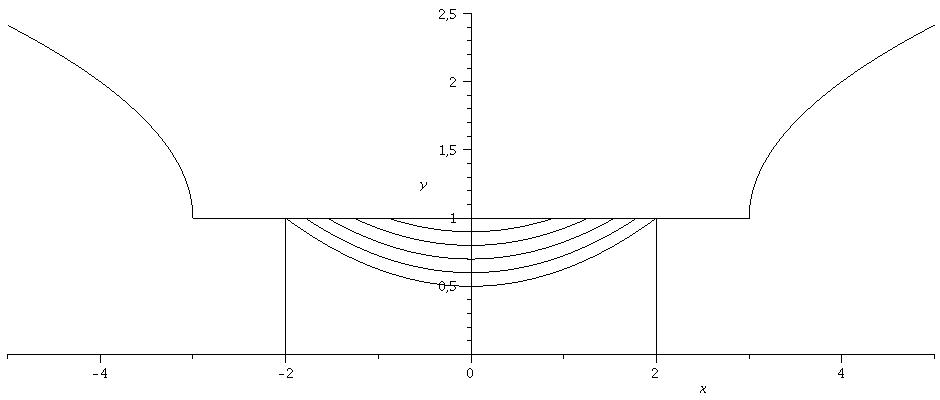}
\end{center}


Note that $\vp_k\ge \vp_{k-1}$. The following remark will be important: at points where $\vp_{k}<1$ we have
$$
\vp_k= \vp_{k-1}+\frac12 4^{-k}.
$$
We also have
$$
\inf \vp_k=\vp_k(0)> 1/2 \quad \mbox{  for } \ k\ge 1.
$$
Moreover, $\vp_\infty(x)\le 1$ precisely for $|x|\le \sqrt{2}$ and $\vp_1(x)< 1$ for $|x|< 2$. This means in particular that $\vp_k(x)=1+(|x|^\ve-2)_+$ for $|x|\ge 2$, $k\ge 1$.
Moreover, $\vp_\infty(x)=(|x|^2+14)/16 \le 15/16$ for $|x|\le 1$.

A more general version of the same construction takes
 \begin{equation}
\overline{\vp}_k= 1 -\frac1{2C} + \frac{|x|^2}{4C}-\frac12 C^{-k}
\end{equation}
with $C$ possibly larger than $4$. In that case $1-\vp_\infty(x)\ge 1/4C$ for $|x|\le 1$.

\medskip

\noindent  For the rest of this proof we write $u_k^+=(u-\vp_k)^+\ge0$, $u_k^-=(u-\vp_k)^-\le0$. Notice that the support of $u_k^+$ is contained in the ball of radius 2 as a consequence of the assumption (\ref{est.apriori1}).

We are ready to tackle the RHS of  Identity \eqref{energ.id} with this choice of test functions. One part will be controlled by a small multiple of the present energy, i.\,e., we will absorb it into the LHS of \eqref{energ.id}. The rest will be bounded above by a large multiple of $|\{u_k^+>0\}|$ (the notation $|.|$ means the measure of the set). We recall our goal: if we do this, together with Sobolev inequality, we will get an iterative relation
for the LHS energies
$$
{\cal A}_{k+1}\le C^k ({\cal A}_k)^{1+\sigma}
$$
that converges to zero as $k\to \infty$, as desired, if the iteration is started from a small initial value $A_0$.

\medskip\newpage

\noindent {\bf \ref{sec51}.4. Estimate of the remaining $\cal B$  term}

We start the process  with ${\cal B}_2={\cal B}(u_k^+,\vp_k).$ By inspecting the integral we easily get
$$
{\cal B}_2 \le \gamma {\cal B}(u_k^+,u_k^+)+\frac1{\gamma}{\cal B}^*(\vp_k,\vp_k)\,,
$$
for every $\gamma>0$, where ${\cal B}^*(\vp_k,\vp_k) $ indicates that the integral is performed only on the set where either $x$ or $y$ belong to  $\{u_k^+>0\}$. That is,
$$
{\cal B}^*=\diint [\chi_{\{u_k^+>0\}}(x) + \chi_{\{u_k^+>0\}}(y)  ]K(x-y)(\vp_k(x) -\vp_k(y))^2.
$$
For $\gamma $ small, then $\gamma {\cal B}(u_k^+,u_k^+)$ is absorbed into the LHS (into the energy). Now, using that
$$
|\varphi_k(x)-\varphi_k(y)| \le C \min (1, |x-y|),
$$
and the size of the kernel $K$, we arrive to the estimate
$$
{\cal B}^*
\le C |\{u_k^+>0\}|\le C 4^{2k}\dint_{\{u^+_k>0\}} (u^+_{k-1})^2dx
$$
 The last inequality follows by Chebyshev's inequality, since $u_{k-1}\ge 4^{-k}/2$ whenever $u^+_k>0$. The obtained expression is good for our later purposes.

\medskip

\noindent {\bf \ref{sec51}.5. Analysis of the ${\cal Q}$ terms for $0<s<1/2$}

The  last term in \eqref{energ.id} also has a bilinear structure. Indeed,
$$
\begin{array}{l}
\diint u_k^+(x)\frac{\nabla \vp_k(x)}{\vp_k(x)} \nabla L (x-y)[u(x)-u(y)]\,dxdy := {\cal Q}(u_k^+,u) \\[14pt]
={\cal Q}(u_k^+,u_k^+) + {\cal Q}(u_k^+,\vp_k)+{\cal Q}(u_k^+, u_k^-)=
{\cal Q}_1+{\cal Q}_2+{\cal Q}_3\,,
\end{array}
$$
but note that the ``kernel" that is involved is not symmetric due to the presence of terms with $\vp_k$. The study of the  contribution of each of the three terms is again split into the
close-range and far-field interactions, represented by the integrals for $x-y$ lying in a ball around the origin, or in its complement. In that sense we note that  $\nabla L$ satisfies $|\nabla L|\le c|x-y|\,K(x,y)$. This will be used repeatedly.

\noindent $\bullet$ Let us first tackle the  integral ${\cal Q}_1$ in a ball of radius $\eta$ around  the origin:
$$
\begin{array}{l}
|{\cal Q}(u_k^+,u_k^+)^{int}|\le \diint_{|x-y|\le \eta}  u_k^+(x)\left|\frac{\nabla \vp_k(x)}{\vp_k(x)}\right| \,c\,|x-y|\,K(x,y)|u_k^+(x)-u_k^+(y)|\,dxdy\\[12pt]
\le 8c^2 \diint_{|x-y|\le \eta}  (u_k^+)^2(x)\left|{\nabla \vp_k}/{\vp_k}\right|^2 |x-y|^2K(x,y)\,dxdy\\+
\frac1{4}\diint_{|x-y|\le \eta}  K(x,y)|u_k^+(x)-u_k^+(y)|^2\,dxdy\,.
\end{array}
$$
The last integral is just the part of $(1/4){\cal B}(u_k^+,u_k^+)$ integrated for $|x-y|\le \eta$, so it can be absorbed by the LHS energy. The first integral is first integrated in $y$ which is easy since
$$\int_{|x-y|\le \eta}K(x,y)|x-y|^2dy=O(\eta^{2s}).
$$
 Using this and also that $|\nabla \vp/\vp| \le C$, this  first integral can be estimated then as
$$
\le C\eta^{2s} \dint (u_k^+)^2\,dx \,,
$$
an expression that can be left in the RHS or absorbed into the LHS if we take $\eta$ small.

   \noindent $\bullet$ Let us compute the outer part of ${\cal Q}_1$ (for $|x-y|>\eta$).   In this region, $\nabla L$ is
integrable, so that
$$
\big |{\cal Q}^{out}(u_k^+,u_k^+)\big | \le C(\eta)  \left(\dint (u_k^+)^2dx+(\dint u_k^+\,dx)^{2}\right)\,,
$$
which is also admissible, as we will see. Note that the last integral comes from the term in $u_k^+(x)u_k^+(y)$.

\noindent $\bullet$ Next, we treat the term
$$
{\cal Q}_2={\cal Q}(u_k^+,\vp_k)=\diint u_k^+ \frac{\nabla \vp_k}{\vp_k} \nabla L (x-y)[\vp_k(x)-\vp_k(y)]dxdy.
$$
Remember that $u_k^+(x)$ is compactly supported in a small ball.  For $|x-y|\le 4$ we have $|\vp_k(x)-\vp_k(y)|\le C|x-y|$, so
$$
\nabla L (x-y)[\vp_k(x)-\vp_k(y)]
$$
is integrable in $y$ and we are left with
$$
C \dint u_k^+(x)dx\le C |\{u_k^+>0\}|\le C^k \dint_{\{u^+_k>0\}}  (u^+_{k-1})^2dx,
$$
which is also a good term.

\noindent $\bullet$  For $|x-y|$ larger the calculation is more involved. We will use the fact that $\nabla L$ has mean value zero on spheres, Therefore, we write \ $\vp_k(x)-\vp_k(y)=(\vp_k(x)-\vp_k(y-x))+ \linebreak (\vp_k(y-x)-\vp_k(y))$. We observe
that
$$
\dint \nabla L(x-y)(\vp_k(x)-\vp_k(y-x))\,dy=\dint_{|z|\ge 4} \nabla L(z)(\vp_k(x)-\vp_k(z))\,dz
$$
is zero in the sense of principal value since (i) $(\vp_k(x)-\vp_k(z))$ is a radial function of $z$, (ii) we have an antisymmetry property for $\nabla L$; both facts imply the cancelation of the integral. The rest of the integral is
$$
\dint \nabla L(x-y)(\vp_k(y-x)-\vp_k(y))\,dy
$$
Since $|\nabla L(x-y)|\sim |x-y|^{-(N+1-2s}$ as $|y|\to\infty$ and $|\vp_k(y-x)-\vp_k(x)|\sim C|x||y|^{\ve-1}$ the whole integral is convergent if $2>2s+\ve$, which is a smallness condition on $\ve_0$.
Performing then the integral in $x$, we get the conclusion that
$$
|{\cal Q}(u_k^+,\vp_k)| \le C |\{u_k^+>0\}|\le C^k \dint_{\{u^+_k>0\}}  (u^+_{k-1})^2dx\,,
$$
as desired.

\noindent $\bullet$ The last term to examine is
\begin{equation}\label{delicate}
{\cal Q}_3={\cal Q}(u_k^+, u_k^-)=-\diint u_k^+(x)\frac{\nabla \vp(x)}{\vp(x)} \nabla L (x-y)u_k^-(y)\,dxdy.
\end{equation}
For $|x-y|\le \eta$ small, we use that $|\nabla L (x-y)|\le C |x-y|\,K (x-y)$ and then $Q(u_k^+,u_k^-)$ is bounded by a small fraction of  ${\cal B}(u_k^+,u_k^-) $ (remember that this term had the right sign). We can therefore get this part absorbed  by the LHS of the energy identity.

\noindent $\bullet$ Finally, for $|x-y|>\eta $, we have the worst convergence case. This is the only place where we use the restriction $s<1/2$. We solve the difficulty of the integrability in $y$ at infinity by taking  $\ve < \ve_0=1-2s$, so that integration first in $y$ is bounded, since the term $|\nabla L(x-y)u_k^-(y) |$ is of the form     $O(|y|^{2s+\ve-N-1})$. We conclude that
\begin{equation}\label{delicate2}
|Q^{out}(u_k^+,u_k^-)|\le C|\{u_k^+>0\}|\le C^k \dint_{\{u^+_k>0\}} (u_{k-1}^+)^2.
\end{equation}

\noindent {\bf Summary}. Using \eqref{ineq.H1}, we  obtain  for $0<s<1/2$ and $t_1<t_2\leq 0$ \ the following energy inequality:
\begin{equation}\label{mod.ener.ineq}
\begin{array}{l}
\dint \frac{(u_k^+ (t_2))^2}{\vp_k}dx + \frac12 \dint_{t_1}^{t_2}{\cal B}(u_k^+,u_k^+)\,dt\\[16pt]
\le
 2\!\dint \frac{(u_k^+(t_1))^2}{\vp_k }dx + C^{2k}\dint_{t_1}^{t_2}\!\!\dint_{\{u^+_k>0\}} (u_{k-1}^+)^2\,dx dt\,,
\end{array}
\end{equation}
where $C$ is a universal constant that only depends  on  $s$ and the dimension, $N$.
In the application to the iteration, the $t_i$ will be chosen in dependence of $k$.

\medskip

\noindent {\bf \ref{sec51}.6. Iteration and end of proof of Lemma \ref{reg.1}}

This part is very similar to the one at the end of the boundedness proof in section \ref{subs.endbound}.
We define the total energy function for the truncated solution $u_k^+$ as
\begin{equation}\label{def.aj}
{\cal A}_k=\sup_{ T_k\le t\le 0}\int  (u_k^+)^2(t)\,dx + \int_{T_k}^0 {\cal B}( u_k^+, u_k^+)\,dt,
\end{equation}
where $T_k=-2(1+2^{-k})$, $k=0,1,\cdots$. Notice that $\vp_k$ lies between $1/2$ and 1 at the points where $u_k^+$ is not zero. From (\ref{mod.ener.ineq}) with $k\ge 1$, taking arbitrary values $t_2=t\ge T_k$ and $t_1=t' \in [T_{k-1},T_k]$ we have
\begin{equation} \label{bdd2}
{\cal A}_k\leq 4\inf_{t' \in [T_{k-1},T_k]} \int  (u_{k}^+)^2(t')\,dx+C^{2k}\dint_{t'}^{0}\!\!\dint_{\{u_k^+>0\}} (u_{k-1}^+)^2\,dx dt=I+II\,.
\end{equation}
Taking averages in $t'$ we arrive at the inequality
$$
\inf_{t' \in [T_{k-1},T_k]} \int  (u_{k}^+)^2(t')\,dx
\leq \frac 1{T_k-T_{k-1}}\int_{T_{k-1}}^{T_k} \int  (u_{k}^+)^2(t')\,dx dt'
$$
$$
\le 2^{k}  \int_{T_{k-1}}^{T_k} \int  (u_{k}^+)^2(t')\,dx dt'.
$$
Observing that $u_{k}^+(x)>0$ implies $u_{k-1}^+(x)>u_{k}^+(x)+4^{-k}/2$, we can realize that both, $I$ and $II$, have the same flavour, and that in fact we have the estimate
\begin{equation} \label{recur}
{\cal A}_k\leq C^{k}\dint_{T_{k-1}}^{0}\!\!\dint_{\{u_{k-1}^+>4^{-k}/2\}} (u_{k-1}^+)^2\,dx dt\,,
\end{equation}
for a possibly larger constant $C$.  The next step it to modify the proof in section \ref{subs.endbound}, replacing the Sobolev inequality by the second part of Lemma \ref{lem.q}. To that end, let $p>2$ be the exponent corresponding to Sobolev's embedding theorem so that
$$
\dint (u_{k-1}^+)^p\,dx \leq C \left[{\cal B}( u_{k-1}^+, u_{k-1}^+)\right]^{p/2}.
$$
Take $\theta=2/p$ and define $q=(1-\theta)2+\theta p$. Then
$$
\begin{array}{ll}
\dint_{\{u_{k-1}^+>4^{-k}/2\}} (u_{k-1}^+)^2\dx \le 4^{(k+1)(q-2)}\dint (u_{k-1}^+)^q\, dx \\
\le 4^{(k+1)(q-2)}
\left(\dint (u_{k-1}^+)^2\,dx\right)^{(1-\theta)}\left(\dint (u_{k-1}^+)^p\,dx\right)^{\theta}
\\
\le C  4^{k(q-2)} \left(\dint (u_{k-1}^+)^2\,dx\right)^{(1-\theta)}{\cal  B}(u_{k-1}^+,u_{k-1}^+)
\end{array}
$$
Integration in time $t$ along the interval $[T_{k-1}, \, 0]$ gives us from inequality \eqref{recur} and the previous estimate a recurrence relation  of the form
$$
\begin{array}{l}
{\cal A}_{k} \le   C^k\left( \sup_{ T_{k-1}\le t\le 0}\dint  (u_{k-1}^+)^2(t)\,dx\right)^{1-\theta} \cdot \dint_{T_{k-1}}^0 {\cal B}( u_{k-1}^+, u_{k-1}^+)\,dt \\
\le  C^k {\cal A}_{k-1}^{(1-\theta)}  {\cal A}_{k-1} =C^k {\cal A}_{k-1}^{1+\tau} ,
\end{array}
$$
with $\tau=1-\theta>0$ and a possibly larger constant $C$.

 We need $\delta$ to be very small to start the iteration so that the sequence ${\cal A}_{k}$ converges and then  ${\cal A}_\infty=0$, which means that $u\le \eta_\infty$ and this in turn implies that  $u\le 7/8$ for $|x|\le 1$. We thus get the result in the Lemma statement with $\mu=1/8$.

\medskip

 \noindent {\bf Remark.} A simple modification of $\vp_\infty$ would give other values of $\mu\in (0,1/2)$, of course with a  different estimate of the maximum allowed value for $\delta$. The proof also shows that the time size $T=4$ can be replaced by any other number and the iteration will work with a different value for $\delta$ (and the same values for $\mu$ and $\ve$). \qed

 \section{Modification of the energy calculation}\label{sec-mod}

In the iterative process that we will consider below there will be situations in which the solution considered in a cylinder as above is bounded between two positive constants $0<M_1<M_2$. We want to establish that a similar result holds and the relation $\delta$-$\mu$ does not change much, which will essential in our iterations. The use of  rescaling and the translation invariance in $(x,t)$ allow to recover a solution defined in the standard domain $\Gamma_4$ which is the one chosen for all our calculations. But imposing the normalization $0\le u\le 1$ asks for a vertical translation in $u$ to adjust the lower level (on top of the usual scaling), and this leads to a modified  equation with the following form
\begin{equation}\label{new.eqn}
\partial_t u =\nabla \cdot(D(u)\nabla {\cal L}u)
\end{equation}
where $D$ has the form $D(u)=d_1 + d_2 u$. We will normalize so that $d_1+d_2=1$ (i.\,e.,  $D(1)=1$) and we will have $|d_2|< 1$ (in practice, $d_2$ becomes small as the iterations advance).

We re-do the energy calculations of the previous section but this time we use as test function
$\eta= F(u\vee \varphi)-F(\vp)$, where $F$ is defined as
$$
F(u):=\int_1^u \frac1{D(s)}\,ds=\frac1{d_2}\log (d_1+d_2u)
$$
Note that $\eta\ge 0$ and $\eta=0$ where $u\le \vp$. Then the LHS gives
$$
I_{lhs}= \dint_{T_1}^{T_2}\dint \eta\,\partial_t u\,dxdt =
\frac1{d_2}\dint_{T_1}^{T_2}\dint \log\left(\frac{d_1+d_2 (u\vee \vp)}{d_1+d_2  \vp}\right)\partial_t u \,dxdt\,.
$$
Next, we write
$$
\frac{d_1+d_2 (u\vee \vp)}{d_1+d_2  \vp}=1+\beta u_\vp^+, \quad \beta=\frac{d_2}{d_1+d_2  \vp}\,,
$$
and note that the function
$$
H(s):=\frac1{d_2}\{s\log(1+\beta s)+\frac1{\beta}\log(1+\beta s)-s\}
$$
is such that $H'(s)=d_2^{-1} \log(1+\beta s)$, so that by integration in time we arrive at
$$
I_{lhs}=\int H(u_\vp^+(x,T_2))\,dx-\int H(u_vp^+(x,T_1))\,dx.
$$
We  will need an estimate of this quantity. Since $H(0)=H'(0)=0$ and $H''(s)=\beta/(d_2(1+\beta s))$, i.\,e., $H''(u_\phi^+)=1/(d_1+d_2 (u\vee \vp))$,  we have  the estimate $H(u_\phi^+)\sim (u_\phi^+)^2$
at all points where $u\ge \vp$, $u$ is bouded above and $\vp$ is bounded below away from zero in the sense that
\begin{equation}\label{ineq.H1.b}
c_1\, (u_\vp^+)^2 \le H(u_vp^+)\le c_2 (u_\vp^+)^2\,.
\end{equation}
and the constants go to 1/2 as $d_2\to 0$ (and $d_1\to 1$), since in the limit $H''(s)=1$.

On the other hand, on the  RHS of the weak formulation, instead of Formula  \eqref{star}, we have the following:
$$
\begin{array}{l}
\dint dt \dint (F(u\vee \varphi)-F(\vp))\,\, \mbox{\rm div} \left[ D(u(x))\, \nabla L (x-y)\, (u(x)-u(y))\right] \,dxdy\\
=-\dint dt \dint_{u>\vp}\{\frac{\nabla u(x)}{D(u(x))}-\frac{\nabla \vp(x)}{D(\vp(x))} \}\, D(u(x)) \nabla_x L (x-y)\, (u(x)-u(y)) \,dxdy\,,
\end{array}
 $$
which we again split as ${\bf I}+{\bf II}$. In the present situation we take
$$
{\bf I}= -\dint dt \diint \nabla u_\vp^+(x) \nabla_x L (x-y)[u(y)-u(x)]\,dxdy.
$$
After integrating by parts and symmetrizing, we get $-{\bf I}= (1/2)\int  {\cal B}(u_\vp^+,u)\,dt$, where
$$
{\cal B}(u_\vp^+,u).=\diint (u_\vp^+(x)-u_\vp^+(y))\, K(x-y) (u(x)-u(y))\,dxdy\,.
$$
As before, we  separate the expression into three integrals, using the splitting: $u=u_\vp^++\vp+u_\vp^-$.
The rest of the integral in the RHS takes the form
$$
{\bf II}= \dint dt \diint_{u(x)>\vp(x)} \frac{D(u)-D(\vp)}{D(\vp)}\nabla \vp(x)\nabla L(x-y) (u(x)-u(y))\,dxdy.
$$

\noindent $\bullet$  We note that when $D(u):=d_1+d_2u$ with $d_1+d_2=1$, \ we have \ $F(s)=(1/d_2)\log(d_1+d_2u)$ and \ $
{\bf II}= \dint  \widehat {\cal Q}(u_\vp^+,u)\,dt,$ with
$$
\widehat {\cal Q}(u_\vp^+,u)= d_2 \diint u^+_\vp(x)\frac{\nabla \vp(x)}{D(\vp(x))}\nabla L(x-y) (u(x)-u(y))\,dxdy\,,
$$
which looks like $\cal Q$ of previous section but for an interesting small factor, $d_2$.

Finally, the energy inequality takes  the form
\begin{equation}\label{energ.ineq}
\left\{\begin{array}{l}
\dint H(u_\vp^+) \,dx \Big|_{T_1}^{T_2}+ C\dint_{T_1}^{T_2}{\cal B}(u_\vp^+,u_\vp^+)\,dt+
C \dint_{T_1}^{T_2}{\cal B}(u_\vp^+,u_\vp^-)\,dt\\[10pt]
\le C_2 \,\dint_{T_1}^{T_2}{\cal B}(u_\vp^+, \varphi)\,dt
+\dint_{T_1}^{T_2} \widehat {\cal Q}(u_\vp^+,u)\,dt\,,
\end{array}\right.
\end{equation}
to be compared with identity \eqref{energ.id}.

 Repeating the rest of the steps of the previous section offers no novelties and we arrive at a similar result

\begin{lem} \label{reg.1.mod} The result of Lemma {\rm \ref{reg.1}} holds for the weak energy solutions of Equation {\rm \eqref{new.eqn}} under the assumptions $D(u)=d_1+d_2u$, $d_1,d_2\ge 0$, $d_1+d_2=1$, $d_1\ge 1/2$, and he result holds with same constants $\mu$, $\delta$ for different values of $d_1, d_2$.
\end{lem}

\section{Pulling up from zero. Proof of Lemma \ref{reg.1b}}\label{sec53}

Here, we are trying to pull the equation uniformly away from zero under some assumption on the size of a level set. The precise assumptions are
$$
0\le u\le \max\{1, |x|^\ve -1\}\quad \mbox{ in \ } \ \ren\times [-4,0],
$$
and  $ |\{u\le 1/2\}\cap \Gamma_4|\le \delta |\Gamma_4|. $ The desired conclusion is then that \ $ \left. u\right|_{\Gamma_2}\ge \mu_0>0$, for an appropriate $\delta>0$ to be chosen.

The tools are energy inequalities and integration by parts. In order to deduce the necessary energy inequalities we recall the bilinear forms and  integration by parts formulas of Section \ref{sec51}. We use the same notation for the kernels $K$ and $L$.

\medskip

\noindent {\bf \ref{sec53}.1. Local energy inequality.}  The basic calculation is as follows. We take a positive smooth function $\varphi$  and use  $\eta=\log(v)$ as a test function in  the weak formulation, where $v(x,t)$ is defined as $v=u/\varphi$ if $u<\varphi$ and $v=1$ otherwise. In other words,
\begin{equation}\label{func.v}
v=\frac{u}{\varphi}\wedge 1, \quad \mbox{or \ } \ v-1=\frac{(u-\varphi)^-}{\varphi}:=\frac{u_{\varphi}^-}{\varphi}\,.
\end{equation}

 Note that $0\le v \le 1$. We will use the notations: $u_\varphi^+=(u-\varphi)^+$, $u_{\varphi}^-=(u-{\varphi})^-$, so that one has \ $ u=\varphi+ u_{\varphi}^+ + u_{\varphi}^-.$ Next, we put
$$
h(v)= \big[ v\log(v)-(v-1)\big], \quad \forall v\ge 0\,,
$$
This function takes as minimum $h(1)=0$ and $h(v) > 0$ for all $v>0, v\ne 1$ with $h(0)=1$. Moreover,
$h''(v)=1/v$ so that for $0<v\le 1$ we have $h''(v) \ge 1$  and
\begin{equation}\label{ineq.h}
\,h(v)= \big[ v\log(v)+1-v\big]\ge \frac{1}{2}  (v-1)^2=\frac{(u_{\varphi}^-)^2}{2\varphi^2}.
\end{equation}
We also have the inequality
\begin{equation}\label{ineq.h2}
\,h(v) \le  (v-1)^2=\frac{(u_{\varphi}^-)^2}{\varphi^2}.
\end{equation}
This is used below. Next, we have the following calculation for the choice of function $v=v(x,t)$ made in \eqref{func.v}
$$
\frac{d}{dt}\int \varphi \,h(v)\,dx=\int \log(v) \varphi \,v_t= \int_{u< \varphi} \log(v) \varphi\, v_t= \int \log(v)u_t\, ,
$$
since $\log(v)=0$ for $u\ge \varphi$. Using the weak formulation of the equation,     we now get
\begin{equation}\label{lem.ee}
\int \varphi \,h(v(t_2))\,dx + \int_{t_1}^{t_2}\int \nabla \big[\log(v)\big]\,u\nabla {\cal L}u\,dxdt=
\int \varphi \,h(v(t_1))\,dx.
\end{equation}
Let us work out the meaning of the middle term (energy term). First, we have
\begin{equation}\label{energy.disp}
\begin{array}{l}
\dint \nabla \big[\log(v)\big]\cdot u\nabla {\cal L}u\,dx=
\dint  u \frac{\nabla v}{v}\,\cdot \nabla {\cal L}u\,dx=\dint_{\{u<\varphi\}}  \varphi \nabla (u_\varphi^-/\varphi)\,\cdot \nabla {\cal L}u\,dx\\ \hskip 4.4cm=
\dint_{\{u<\varphi\}}  \nabla u_\varphi^-\,\cdot \nabla {\cal L}u\,dx-
\dint_{\{u<\varphi\}}  \frac{u_\varphi^-}{\varphi} \nabla \varphi\,\cdot \nabla {\cal L}u\,dx\,.
\end{array}
\end{equation}
Clearly, using $u=\varphi+ u_{\varphi}^+ + u_{\varphi}^- = (u\vee \varphi) + u_{\varphi}^-$, we get
$$
\dint_{\{u<\varphi\}}  \nabla u_\varphi^-\,\cdot \nabla {\cal L}u\,dx={\cal B}(u_\varphi^-,u)=
{\cal B}(u_\varphi^-,u_\varphi^-)+{\cal B}(u_\varphi^-,u_\varphi^+)+{\cal B}(u_\varphi^-,\varphi).
$$
The first term is the one we want to keep to complete the expression of the energy.
The second term  ${\cal B}(u_\varphi^-,u_\varphi^+)$ is also positive in view of formula \eqref{eq.bilin} since both functions have opposite signs and disjoint supports, cf. Corollary \ref{cor3.3}.   Hence, this term could be discarded, but it will turn out to be useful as we have already seen.

The remaining term ${\cal B}(u_\varphi^-,\varphi)$ and the last integral in \eqref{energy.disp} are delicate
and it is there that we have to make an  argument with a careful choice of test function. Summing up, we have
\begin{equation*}
\dint \nabla \big[\log(v)\big]\cdot u\nabla {\cal L}u\,dx \ge
{\cal B}(u_{\varphi}^-, u_{\varphi}^-)+ {\cal B}(u_\varphi^-,\varphi)+ {\cal B}(u_\varphi^-,u_\varphi^+)
-\dint_{u<\varphi}  \frac{u_\varphi^-}{\varphi} \nabla \varphi\,\cdot \nabla {\cal L}u\,dx\,.
\end{equation*}
Putting
$$
 {\cal Q}(w,u)=\diint w(x)\frac{\nabla \vp(x)}{\vp(x)} \nabla L(x-y)[u(y)-u(x)]dxdy,
 $$
we arrive from the above and from \eqref{ineq.h} and \eqref{ineq.h2} to the  basic energy inequality:
\begin{equation} \label{ener}
\begin{array}{c}
 \dint \frac{1}{2\vp}\,(u_\vp^-(t_2))^2 dx + \dint_{t_1}^{t_2}{\cal B}(u_\vp^-,u_\vp^-)\,dt+
\dint_{t_1}^{t_2}{\cal B}(u_\vp^-,u_\vp)\,dt\\[10pt] \le \dint \frac{1}{\vp}\,(u_\vp^-(t_1))^2 dx-\dint_{t_1}^{t_2}{\cal B}(u_\vp^-,\varphi)\,dt
+\int_{t_1}^{t_2}  {\cal Q}(u_\vp^-,u)\,dt.
\end{array}
\end{equation}
The local energy,  in the time interval $[T_1,T_2]$,  is now defined as
\begin{equation}
E_\varphi(u):= \sup_{T_1<t<T_2}\int \frac1{\varphi}[u_\varphi^-]^2\,dx + \int_{T_1}^{T_2} {\cal B}(u_\varphi^-,u_\varphi^-)\,dt.
\end{equation}
Note that the test function $\log v$ is negative but we are interested in $u$ being a supersolution, and this is the case for instance if we truncate it by 2: $\bar u=u\wedge 2$. So, we do not need to worry about growth at infinity.

\medskip


{\bf \ref{sec53}.2. The iterative process}

\noindent {\bf (i)} At this point, following De Giorgi's idea we want to obtain an iterative relation playing the energy inequality against the Sobolev embedding using a convenient sequence of cutoff functions. We consider a series of cutoffs with dyadic separation, so that
$$
\varphi_{k+1} \le \varphi_k -a 2^{-k}
$$
holds in the support of $\varphi_{k+1}$ for some fixed $0<a<1$. Moreover, the $\varphi_k$ converge to $\mu_0$  from above in $B_2(0)$. We will use as a test function  $\eta=\log v$ with  $v=u/\varphi_k$ if $u<\varphi_k$ and $v=1$ otherwise. We arrive at the energy inequality at the end of last subsections with $\varphi_k$ instead of $\varphi$.

\medskip

\noindent {\bf (ii)} We propose  a concrete construction of the cutoffs $\varphi_k$. All functions $\varphi_k(x)$ are symmetric, nonnegative, non-increasing and compactly supported and the sequence
is nonincreasing with $k$. The basic profile $\varphi_0$ is a kind of rounded mesa:
$$\varphi_0(|x|\left\{
\begin{array}{l}
=\frac 12[4-|x|]^m \quad \quad   \mbox{when \ }7/2 <|x| < 4,\\
\equiv 1 \qquad \qquad \qquad  \mbox{for \ } \ |x| \le 3,\\
\mbox{a ${\cal C}^2$ decreasing radial function for $3 \le |x|<4$,}
\end{array}\right.
$$
for some $m\ge2$ to be chosen later.

- To construct $\varphi_k$ we first rescale the graph of $\varphi_0$ from the interval $[3,4]=[2+1,2+2]$ to the interval $[2+2^{-k}, 2+2^{-k+1}]$, and extend by $\mu_0+ (1- \mu_0)2^{-k}$ inside the ball $B_{2+2^{-k}}$.
Then, $\varphi_k$ has the following properties:
$$
\begin{array}{ll}
(a) \qquad \qquad \qquad \qquad &\vp_k(x)\le \vp_{k-1}(x)\le \cdots \vp_{0}(x)\,,\\[1pt]
(b)  \qquad \qquad \qquad \qquad &|\nabla \vp_k|/\vp_k\le C^{k}\varphi_k^{-1/m}\,, \\[1pt]
(c) & \varphi_{k-1}-\varphi_{k}\ge (1-\mu_0)2^{-k} ,\quad \mbox{in the support of $\varphi_k$},\\[1pt]
(d) & \vp_k \to \mu_0\chi_{B_2} \quad \mbox{ as } k \to \infty,
\end{array}
$$
so that $a=1-\mu_0$.
\begin{center}
\includegraphics[width=12 cm]{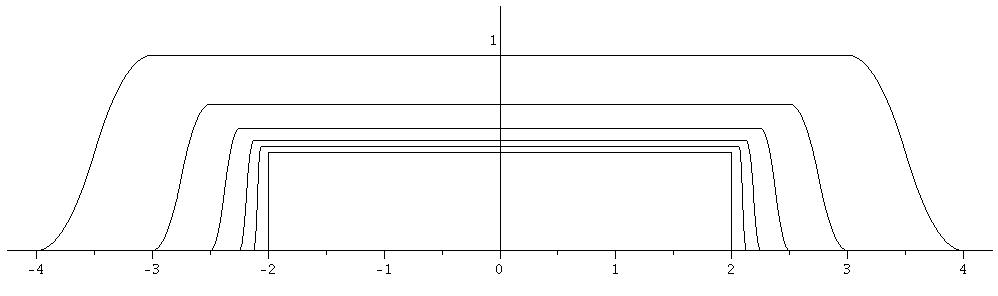}
\end{center}
An explicit expression for $\varphi_k$  could be the following one:
$$
\varphi_k(x)=(\mu_0+ (1- \mu_0)2^{-k})\cdot\vp_0(2^k|x|-2^{k+1} +2).
$$


\

\noindent {\bf (iii) Sobolev embedding}. The embedding that we need is a variant of the one used in preceding sections. We have seen that ${\cal B}(u^-_\varphi,u^-_\varphi)$ is the square $H^r$ norm of $u^-_\varphi$, and in this way it controls an $L^p$ norm of $u^-_\varphi$ for a given $p>1$ that we have calculated.
Coupling this with the energy term
$$
\sup_{T_1<t<T_2}\int \frac1{\varphi}[u^-_\varphi]^2\,dx\,,
$$
we can control for some $q>2$ and $\theta>0$
$$
\big(\int_{T_1}^{T_2}\int \varphi^{-\theta}|u^-_\varphi |^q\,dxdt\big)^{2/q}.
$$
(Since the cutoff $\varphi$ goes to zero, the term $\varphi^{-\theta}$ will be important in controlling the term containing $\nabla\varphi/\varphi$). The proof is as follows:
Write $q$ as a convex combination, $q= 2\theta +p\, (1-\theta)$. Then,
$$
\int \varphi^{-\theta}|u^-_\varphi |^q\,dx\le
\big(\int \varphi^{-1}|u^-_\varphi |^2\,dx\big)^\theta\big(\int |u^-_\varphi |^p\,dx\big)^{(1-\theta)}
$$
We make the choice  $\theta= (p-2)/p$ with $p$ as the $q$ in Lemma \ref{lem.q}. Then, $1+\theta=q/2$ and
$$
\left(\int_{T_1}^{T_2}\int \varphi^{-\theta}|u^-_\varphi |^q\,dxdt\right)^{2/q}\le E_\varphi.
$$
From now and on we will write $u_k^{+/-}$ and $E_k$ to denote  $u_{\varphi_k}^{+/-}$ and $E_{\varphi_k}$ respectively.

\

\noindent {\bf \ref{sec53}.3. Iterated energies.}
In the left of \eqref{ener} we truncate in time along an increasing sequence $T_k\to -2$ and get, for every $T_{k-1}<t_1<T_k<t_2< 0$,
\begin{equation}  \label{ener2}
\begin{array}{c}
 \dint \frac{1}{2\vp_k}\,(u_k^-(t_2))^2 dx+ \dint_{t_1}^{t_{2}} {\cal B}(u_k^-,u_k^-)\,dt+
 \dint_{t_1}^{t_{2}} {\cal B}(u_k^-,u_k^+)\,dt\\[10pt] \le  \dint \frac{1}{\vp_k}\,(u_k^-(t_1))^2 dx - \dint_{t_1}^{t_{2}} {\cal B}(u_k^-,\varphi_k)\,dt
- \dint_{t_1}^{t_{2}}  {\cal Q}(u_k^-,u)\,dt.
\end{array}
\end{equation}

\

\medskip

{\bf \ref{sec53}.4.  Analysis of the RHS}

We now examine the terms left in the RHS of the energy inequality
\eqref{ener2} that we will call $I, II,$ and $III$ for convenience of reference.
Our purpose is to show that these extra terms are either bounded by a small multiple of the LHS, so that it will be absorbed by the LHS, or they are controlled by a term of the form
$$
C^k\,\int_{T_1}^{T_2}\int \varphi^{-\theta}|u^-_\varphi |^2dxdt.
$$

\noindent {\bf  Estimating {\sl I.}}   As initial times $T_k$ we will choose the time $T_k = -(2+2^{-k})$ and $t_1\in [T_{k-1}, T_k]$ is the point where the value of
$$
 \inf _{T_{k-1}\le t \le T_{k}} \dint \frac{1}{\vp_k}\,(u_k^-(t))^2 dx
$$
is attained. In this way we have
$$
I\le 2^k\int_{-2-2^{1-k}}^{-2-2^{-k}}\int\frac{[(u-\varphi_k)^-]^2}{\varphi_k}\,dx\, dt.
$$
Since $|(u-\varphi_k)^-/\varphi_k|<1$ we have $[(u-\varphi_k)^-]^2/\varphi_k\le |(u-\varphi_k)^-| $, and
we can bound $I$ by
$$
C^k\dint_{T_{k-1}}^{0}\dint  \chi_{\{u_k^-<0\}}.
$$

\noindent {\bf Estimate of {\sl II.}} About this  term we have
$$
\big |{\cal B}(u^-_k,\varphi_k)\big | \le \frac12 \big[{\cal B}(u^-_k,u^-_k) +{\cal B}^*(\varphi_k,\varphi_k)\big]
$$
where the star indicates that we can restrict    the domain of integration defining   ${\cal B}(\varphi_k,\varphi_k)$ to the points   where $u^-_k(x)\ne 0$ or $u^-_k(y)\ne 0$.   So this last term can be replaced by the better expression
$$
\frac12\iint\big[\varphi_k(x)-\varphi_k(y)\big]^2K(x,y)\chi(x,y)\,dxdy,
$$
where $\chi$ is the characteristic function of the set of points $(x,y)$ where either $x\in \supp \ u_{k}^-$ or
$y\in \supp \ u_{k}^-$. From the Lipschitz character of $\varphi_k$     and that $|\nabla\varphi_k|\leq C2^k$    we have for fixed $x$
$$
\int\big[\varphi_k(x)-\varphi_k(y)\big]^2K(x,y)\,dy\le C^k.
$$
Therefore, the above term may be bounded by
$$
C^k \int \chi_{\{u^-_{k}<0\}} dx.
$$
Hence,
$$
II \le \frac 12  \dint_{t_1}^{t_{2}} {\cal B}(u_k^-,u_k^-) dt+ C^k \int_{t_1}^{t_2} \int \chi_{\{u^-_{k}<0\}}\,dxdt.
$$
The first term is absorbed by the local energy in the LHS.

\medskip

\noindent {\bf Estimate of {\sl III.}} The remaining part of the proof is devoted to
the estimate of the term
$$
 {\cal Q}(u_k^-,u)=\diint u_k^-(x)\, \frac{\nabla \varphi_k(x)}{\varphi_k(x)} \,\cdot \nabla L (x-y) \,(u(y)-u(x))\,dxdy
$$
Write, by the usual splitting of $u$,
$$
 {\cal Q}(u_k^-,u)= {\cal Q}(u_k^-,u_k^-)+ {\cal Q}(u_k^-,u_k^+)+ {\cal Q}(u_k^-,\vp_k)={\cal Q}_1+{\cal Q}_2+{\cal Q}_3.
$$

\noindent {\bf ${\cal Q}_1$ :} Noting  that
$$
|\nabla L(x,y)|\le |x-y| K(x,y) =O(|x-y|^{-(N+1-2s)}),
$$
we have from H\"older's
$$
 \left |{\cal Q}_1\right |=\left |{\cal Q}(u_k^-,u_k^-)\right |=\left |\diint u_k^-(x)\frac{\nabla \vp_k(x)}{\vp_k(x)} \nabla L(x-y)[u_k^-(y)-u_k^-(x)]dxdy \right | \leq
 $$
 $$
 \leq \left[\diint \left|u_k^-(x)\frac{\nabla \vp_k(x)}{\vp_k(x)}\right|^2 \frac 1{|x-y|^{N-2s}}dxdy\right]^{1/2}\cdot\left[{\cal B}(u_k^-,u_k^-)\right]^{1/2} \leq
 $$
 $$
 \leq 4 {C_N} \dint \vp_k^{-2/m}(u_k^-)^2 \,dx+\frac 14  {\cal B}(u_k^-,u_k^-),
 $$
where we have used here the properties of the functions $\vp_k$.
 Integration in time gives that the  second term is absorbed by the local energy expression in the LHS and the first is admissible for our purposes.

\medskip

\noindent {\bf Final step.} To study $ {\cal Q}_2$ and $ {\cal Q}_3$ we consider a smooth decomposition of the kernel
 $$
 \nabla L=\nabla L_\psi + \nabla L_{1-\psi},
 $$
 where $1-\psi$ is a bump function supported around the origin. We get
$$
{\cal Q}_2+{\cal Q}_3= \diint u_k^-(x)\frac{\nabla \vp_k(x)}{\vp_k(x)} \nabla L_{1-\psi}(x-y)[u_k^+(y)+\vp_k(y)]dxdy+
$$
$$
+\diint u_k^-(x)\frac{\nabla \vp_k(x)}{\vp_k(x)} \nabla L_\psi(x-y)[u(y)]dxdy =:\bar{\cal Q}_2+\bar{\cal Q}_3.
$$
In other words, in the compact part of the support of $1-\psi$ we keep the expression of $u$ as $u_k^-+\vp_k+u_k^+$, while outside, where $u_k^-=0$ and $\vp_k$ grows, we just keep the term
$$
\dint \nabla L_\psi (x-y) u(y) dy.
$$
\noindent $\bullet$  This term has to  be controlled in a different way (through the change of coordinates if $s\ge 1/2$).  As a consequence, the integrations in $y$ for
 $\bar{\cal Q}_2$ and $\bar{\cal Q}_3$ are convergent and we are left with estimating

$$
\int \left|u_k^-(x)\frac{\nabla \vp_k(x)}{\vp_k(x)} \right| dx  \le   \int \vp_k^{-2/m}(u_k^-)^2 dx +  \int \chi_{\{u^-_{k}<0\}} dx.
$$
The integration in time produces again two admissible terms.

\medskip

\noindent $\bullet$ Finally, for the term
$$
\diint u_k^-(x)\frac{\nabla \vp_k(x)}{\vp_k(x)} \nabla L_{1-\psi}(x-y)[u_k^+(y)]dxdy,
$$
we use the good term
$$
 {\cal B}(u_k^-,u_k^+)= \diint u_k^-(x)K(x-y)[u_k^+(y)]dxdy,
$$
left in the energy inequality, to absorb the integral whenever
$$
\left|\frac{\nabla \vp_k(x)}{\vp_k(x)}\right| |x-y| \le \eta.
$$
In the complement, that is when
$$
\left|\frac{\nabla \vp_k(x)}{\vp_k(x)}\right| |x-y| \ge \eta.
$$
we use that $\vp_k^{-1/m} \ge |\nabla \vp_k|/\vp_k$ and integrate in $y$:
$$
\left | \dint \nabla L_{1-\psi}(x-y)[u_k^+(y)]dy \right | \lesssim \dint_{\eta\vp_k^{\-1/m}}^4 \frac {r^{n-1}}{r^{n+s-1}} dr\sim \max \left(C, \eta \vp_k^{(1-s)/m}\right).
$$
In any case, we end up in the worst of the cases with an expression of the form
$$
\int \vp_k^{-s/m} \left |u_k^-\right | dx,
$$
that we control as before.

\medskip

Inserting  estimates {\sl II} and {\sl III} in \eqref{ener2} we get that
\begin{equation}  \label{ener3}
\begin{array}{c}
 \dint \frac{1}{2\vp_k}\,(u_k^-(t_2))^2 dx+ \dint_{t_1}^{t_{2}} {\cal B}(u_k^-,u_k^-)\,dt
\\[10pt] \le  \dint \frac{1}{\vp_k}\,(u_k^-(t_1))^2 dx + C^k \left[\dint_{t_1}^{t_{2}}\int  \vp_k^{-2/m}(u_k^-)^2 dx dt +  \dint_{t_1}^{t_{2}}\int  \chi_{\{u^-_{k}<0\}} dx dt\right].
\end{array}
\end{equation}
whenever $T_{k-1}<t_1<T_k<t_2<0$. In particular,
\begin{equation}  \label{ener4}
\begin{array}{c}
\sup_{T_k< t_2<0} \dint \frac{1}{\vp_k}\,(u_k^-(t_2))^2 dx+ \dint_{T_k}^{0} {\cal B}(u_k^-,u_k^-)\,dt \le \inf_{T_{k-1}<t_1<T_k}   \dint \frac{2}{\vp_k}\,(u_k^-(t_1))^2 dx
\\[15pt]  + C^k \left[\dint_{T_{k-1}}^{0}\int  \vp_k^{-2/m}(u_k^-)^2 dx dt +  \dint_{T_{k-1}}^{0}\int  \chi_{\{u^-_{k}<0\}} dx dt\right]
\\[15pt]  \le C^k \left[\dint_{T_{k-1}}^{0}\int  \vp_k^{-2/m}(u_k^-)^2 dx dt +  \dint_{T_{k-1}}^{0}\int  \chi_{\{u^-_{k}<0\}} dx dt\right],
\end{array}
\end{equation}
where the last inequality comes from estimate {\sl I}. In fact, the same argument given for this estimate shows that if $m\ge 2$ then the leading term in the last expression is the second one. Thus, if we define
\begin{equation}
{\cal A}_k:=\sup_{T_k< t < 0}\dint (u_k^-)^2 dx+ \dint_{T_k}^0 {\cal B}(u_k^-,u_k^-)\,dt,
\end{equation}
then,  \eqref{ener4} gives

\begin{equation}
{\cal A}_k \le
 C^k  \dint_{T_{k-1}}^{0}\int  \chi_{\{u^-_{k}<0\}} dx dt.
 \end{equation}
 Using that $\vp_{k-1} \ge \vp_k +a2^{-k}$ in the support of $\vp_k$ we get than $|u_{k-1}^-|\ge |u_{k}^-| +a2^{-k}$ in
 ${\{u^-_{k}<0\}} $. Therefore,
 \begin{equation}
{\cal A}_k \le
 C^{(1+q)k}  \dint_{T_{k-1}}^{0}\int |u^-_{k-1}|^q dx dt.
 \end{equation}
The rest follows the same argument as in Lemma \ref{reg.1}: let $p>2$ be the exponent corresponding to Sobolev's embedding theorem so that
$$
\dint |u_{k-1}^-|^p\,dx \leq C \left[{\cal B}( u_{k-1}^-, u_{k-1}^-)\right]^{p/2}.
$$
Take $\theta=(p-2)/p$ and define $q=\theta2+(1-\theta) p$. Then,
\begin{equation}
\begin{array}{ll}
{\cal A}_k \le
 C^{(1+q)k}  \dint_{T_{k-1}}^{0}\left(\int |u^-_{k-1}|^2 dx\right)^\theta {\cal B}( u_{k-1}^-, u_{k-1}^-)dt  \\
  \le   {\widetilde C}^k\left( \sup_{ T_{k-1}\le t\le 0}\dint  (u_{k-1}^-)^2(t)\,dx\right)^{\theta} \cdot \dint_{T_{k-1}}^0 {\cal B}( u_{k-1}^-, u_{k-1}^-)\,dt\\
\le { \widetilde C}^k {\cal A}_{k-1}^{\theta}  {\cal A}_{k-1} ={\widetilde C}^k {\cal A}_{k-1}^{1+\theta} ,
\end{array}
 \end{equation}

This completes the proof of Lemma \ref{reg.1b} for the critical case since, by hypothesis,
$$
A_0=\int_{-4}^0 \int \frac{[(u-\varphi_0)^-]^2}{\varphi_0}\,dxdt\le \iint_{\Gamma_4} \chi_{\{u-\varphi_0<0\}}\le\left| \{u\le 1\}\cap \Gamma_4\right| \le \delta|\Gamma_4|,
$$
and $\delta$ can be chosen as small as we need.

\medskip

\noindent {\bf The noncritical case of Lemma \ref{reg.1b}.} This is the case where the solution lies between, say, $M$ and $M+1$ for some $M>0$. The proof is similar except that  for the lower estimate we may already assume that $u|_{\Gamma_4}\ge \mu>0$ and then all the involved cutoffs  satisfy $\varphi_k>\mu$ so that $\nabla\varphi_k/\varphi_k$ is a smooth bounded function.


\section{\bf The lemma on intermediate values}
\label{sec.intermediate}

 We still need a main ingredient before we attack the regularity issue by means of a suitable iteration. Indeed,  we have to improve Lemma \ref{reg.1} by showing that, in order to get a uniform reduction of the maximum in a smaller ball it is not necessary to ask that $u\le 1/2$ ``most of the time'' in $\Gamma_4$,  but only ``some of the time''. This is precisely stated in Lemma \ref{reg.2},

The proof of this result uses De Giorgi's idea of loss of mass at intermediate levels, which he applied in an elliptic context. In the present  parabolic setting we will follow the idea of the proof of the similar result  that was carried out in the papers \cite{CVass}, \cite{CChV} in the linear framework. We give next the detailed statement of the ''intermediate lemma'' and its proof. We start by selecting a cutoff function of the form $\vp=1+\psi_\lambda+F$, where we choose
$$
F(x)=\sup\{-1,\inf\{0,|x|^2-9\}\}
$$
Note that $F\le 0$ is Lipschitz, compactly supported in $B_3(0)$ and $F=-1$ in $B_2(0)$. Moreover, for $0<\lambda<1/3$ we define
$$
\psi_\lambda(x)=((|x|-\lambda^{-1/\nu})^{\nu}-1)_+ \ \mbox{for } \ |x|\ge \lambda^{-1/\nu}\,,
$$
and $\psi_\lambda(x)=0$ otherwise. This represents a ``wing'' that starts far away when  $\lambda$ is small, as will be the case. A convenient value of $\nu\in (0,1)$ is needed and it will be determined later.  We also define
$$
\vp_1=1+\psi_\lambda+\lambda F, \quad \vp_2=1+\psi_\lambda+\lambda^2 F,
$$
and put $\vp_0=\vp$, so that $\vp_0\le \vp_1\le \vp_2\le 1$ in the ball $B_4(0)$.

\begin{lem}\label{lem.3layer} There exist small constants $\rho, \gamma>0$ depending only on $N$ and $s$, and $\lambda_0\in (0,1)$, depending only on $N$, $s$, $\rho$, $\delta$ such that for any solution $u$ defined in $S_4$ with
$$
u(x,t)\le 1+\psi_\lambda \ \mbox{in } \ S_4,
$$
with $\lambda\le \lambda_0$, and also such that $|\{u<\vp_0\}\cap ((-4,-2)\times B_1)|\ge \rho$,
then we have the following implication: if
$$
|\{u>\vp_2\}\cap ((-2,0)\times \re^N)|\ge \delta\,,
$$
then
\begin{equation}
|\{\vp_0\le u\le\vp_2\}\cap ((-3,0)\times \re^N)|\ge \gamma.
\end{equation}
\end{lem}

The last line asserts that under the stated assumptions the measure of the intermediate level cannot be small

\noindent {\sl Proof.} (i) In our context we start from the energy estimates we have obtained during the proof of Lemma \ref{reg.1}. We have to arrive at a ``correct form'' of the energy inequality.
We recover the energy calculation \eqref{energ.id} with $\vp$ equal to the ``intermediate cutoff'' $\vp_1$, i.\,e.,
\begin{equation}\label{energ.id2}
\dint \vp_1 H(g) \,dx \Big|_{T_1}^{T_2}  +
\dfrac12 \dint_{T_1}^{T_2} {\cal B}((u-\vp_1)_+,(u-\vp_1)_+)dt+
\end{equation}
\begin{equation*}
+ \dfrac12 \dint_{T_1}^{T_2} {\cal B}((u-\vp_1)_+,(u-\vp_1)_-)dt
\le -\dfrac12 \dint_{T_1}^{T_2} {\cal B}((u-\vp_1)_+,\vp_1)dt + \dint_{T_1}^{T_2}{\cal Q}\,dt\,,
\end{equation*}
with $g=1+(u-\vp_1)_+/{\varphi_1}$ and $ H(g)=g\,\log g -(g-1)$. We want to get now estimates on the RHS that show that all the terms are either absorbed by the
LHS or  can be estimated above by $C\lambda^2$, where $C$ is a fixed constant.
With this we will arrive finally to the expression
\begin{equation}\label{Lipschitz}
\begin{array}{l}
\dint \vp_1 (g\,\log g+1 -g) \,dx \Big|_{T_1}^{T_2}  +
c \dint_{T_1}^{T_2} {\cal B}((u-\vp_1)_+,(u-\vp_1)_+)\,dt\\
+ c \dint_{T_1}^{T_2} {\cal B}((u-\vp_1)_+,(u-\vp_1)_-)\,dt
\le  C\lambda^2(T_2-T_1)\,.
\end{array}
\end{equation}
This would complete the preparatory step of the lemma.

\medskip

\noindent (ii) The verification that \eqref{Lipschitz} holds is as follows. For simplicity we write $(u-\vp_1)_+=u^+_1$ and $(u-\vp_1)_-=u^-_1$. Repeating a bit some arguments, to prove \eqref{Lipschitz} we observe  that
$$|{\cal B}(u_1^+,\vp_1)|
 \le \gamma {\cal B}(u_1^+,u_1^+)+\frac2{\gamma}{\cal B}^*(\vp_1,\vp_1)\,,
$$
for every $\gamma>0$, where
$$
{\cal B}^*(\vp_1,\vp_1)=\diint \chi_{B_3}(x)  K(x-y)(\vp_1(x) -\vp_1(y))^2dy\, dx.
$$
For $\gamma $ small, then $\gamma {\cal B}(u_1^+,u_1^+)$ is absorbed into the LHS of  \eqref{energ.id2}. For the second term we have
$$
{\cal B}^*(\vp_1,\vp_1)\leq 4(\lambda^2{\cal B}^*(F,F)+{\cal B}^*(\psi_\lambda,\psi_\lambda))
$$
Since the function $F$ is Lipschitz, compactly supported  and \ $|x-y|^2K(x-y|$ is locally integrable, we have ${\cal B}^*(F,F)\le C$. Now, for $\lambda$ small one has $\lambda^{-1/\nu} \ge 4$ so that $\psi_\lambda(x)=0$ if $|x|<3$. Therefore, using that
$\psi_\lambda(y) \leq |y|^{\nu}$ always, we get
$$
\begin{array}{l}
{\cal B}^*(\psi_\lambda,\psi_\lambda) =\dint_{B_3} \dint_{|y|>\lambda^{-1/\nu}}  K(x-y)(\psi_\lambda(y))^2dy\,dx\\[6mm]
\leq \dint_{B_3} \dint_{|y|>\lambda^{-1/\nu}} \frac {|y|^{2\nu}}{|x-y|^{N+2-2s}} dy\, dx \le C
\dint_{B_3} \dint_{|y|>\lambda^{-1/\nu}} \frac {|y|^{2\nu}}{|y|^{N+2-2s}} dy\, dx \sim \lambda^{\frac {2-2s-2\nu}\nu}.
\end{array}
$$
It suffices to take  $\frac {2-2s-2\nu}\nu\geq 2$, that is, $\nu \le \frac{2-2s}4$, in order to conclude that
$$
{\cal B}^*(\vp_1,\vp_1)\leq {\cal O}(\lambda^2).
$$
We now consider the terms involving $\cal Q$:
$$
{\cal Q}(u_1^+,u)
={\cal Q}(u_1^+,u_1^+) + {\cal Q}(u_1^+,\vp_1)+{\cal Q}(u_1^+, u_1^-).
$$
The main ingredients are the following estimates:

${\bullet} \hskip 1cm  u^+_1 \le \lambda |F|$, since  $0\le u\le 1+\psi_\lambda$. In particular,

${\bullet} \hskip 1cm  u^+_1 \le \lambda \chi_{B_3}$

${\bullet} \hskip 1cm |\nabla L(x-y)| \sim |x-y|K(x-y) \sim \frac 1{|x-y|^{N+1-2s}}$,  integrable at infinity for $s<1/2$.

${\bullet} \hskip 1cm |\nabla \varphi_1(x)|/\varphi_1(x) \leq C$, $ \forall x $, with $C$ independent of $\lambda$, and

${\bullet} \hskip 1cm |\nabla \varphi_1(x)|/\varphi_1(x) \leq C \lambda$, in  the support of $u^+_1$.

\noindent Hence,
$$
\begin{array}{l}
|{\cal Q}(u_1^+,u_1^+) | =\left |  \diint  u_1^+(x)\frac{\nabla \vp_1(x)}{\vp_1(x)} \nabla L(x-y) (u_1^+(x)-u_1^+(y))\,dxdy\right|\\[12pt]
\le C \eta  {\cal B}(u_1^+,u_1^+) +C \lambda \diint_{|x-y|\ge \eta}  \left[u_1^+(x)^2+u_1^+(x)u_1^+(y)\right]|\nabla L(x-y)| \,dxdy\\\le  C \eta  {\cal B}(u_1^+,u_1^+) +C \lambda \dint u_1^+(x)^2dx+ C \lambda \left(\dint u_1^+(x) \, dx\right)^2 \,.
\end{array}
$$
The first term in the last line of the formula goes to the LHS of \eqref{energ.id2} for $\eta$ small and the other two are just
${\cal O}(\lambda^3)$. Also,
$$
{\cal Q}(u^+_1,\vp_1) =\lambda{\cal Q}(u^+_1,F)+{\cal Q}(u^+_1,\psi_\lambda)
$$
Since
$$
|{\cal Q}(u^+_1,F)| \leq  C \lambda \dint  u_1^+(x) \dint |\nabla L(x-y) (F(x)-F(y))|\,dy\,dx
\le C \lambda \dint  u_1^+(x)\, dx,
$$
the first term above is easily seen to be of order ${\cal O}(\lambda^3)$. For the second we have
$$
\begin{aligned}
|{\cal Q}(u^+_1,\psi_\lambda)| & \leq   C \lambda \dint  u_1^+(x) \dint \left |\nabla L(x-y) (\psi_\lambda(x)-\psi_\lambda(y))\right |\,dy\,dx \\[12pt]
 & \leq  C \lambda \dint  u_1^+(x)\dint_{|y|>\lambda^{-1/\nu}} \frac {|y|^{\nu}}{|y|^{N+1-2s}} \,dy\,dx,
\end{aligned}
$$
where we have used that $\psi_\lambda(x)=0$ in the suport of $u^+_1$ and $\psi_\lambda(y)\leq |y|^\nu$. From this we get
$$
|{\cal Q}(u^+_1,\psi_\lambda)|  \leq   C \lambda^{1+\frac {1-2s-\nu}\nu}  \dint  u_1^+(x)\, dx .
$$
Observing that $\nu <1-2s$ from the previous   condition on $\nu$ , we conclude that
$$
{\cal Q}(u^+_1,\psi_\lambda) \leq {\cal O}(\lambda^2).
$$

 Finally,
 $$
\begin{array}{l}
|{\cal Q}(u_1^+,u_1^-) | =\left |  \diint  u_1^+(x)\frac{\nabla \vp_1(x)}{\vp_1(x)} \nabla L(x-y) (u_1^-(x)-u_1^-(y))\,dxdy\right|\\[12pt]
\le C \eta  {\cal B}(u_1^+,u_1^-) + \diint_{|x-y|\ge \eta}   u_1^+(x)\left |\frac{\nabla \vp_1(x)}{\vp_1(x)}\right| |\nabla L(x-y) u_1^-(y)|\,dxdy.
\end{array}
$$
The term $C \eta  {\cal B}(u_1^+,u_1^-)$ is absorbed into the LHS of \eqref{energ.id2} for $\eta$ small. For the other we use the fact that $|u^-_1|\le \varphi_1=1+\lambda F+\psi_\lambda$. Now,
$$
\begin{array}{l}
 \dint  u_1^+(x)\dint_{|y|>\lambda^{-1/\nu}}
 \left |\frac{\nabla \vp_1(x)}{\vp_1(x)}\right| |\nabla L(x-y)(1+\psi_\lambda(y))|\,dxdy \\[12pt]
 \le  C \lambda  \dint  u_1^+(x)\dint_{|y|>\lambda^{-1/\nu}} \frac {1+|y|^{\nu}}{|y|^{N+1-2s}} \,dy\,dx  \leq   C \lambda^{1+\frac {1-2s-\nu}\nu}  \dint  u_1^+(x)\, dx ,
 \end{array}
$$
whereas,
$$
\begin{array}{l}
 \dint  u_1^+(x)\dint_{\{y: |x-y|\ge \eta \;\& \;|y|\le\lambda^{-1/\nu}\}}\left |\frac{\nabla \vp_1(x)}{\vp_1(x)}\right| |\nabla L(x-y)(1+F(y))|\,dxdy \\[12pt]
 \le  C \lambda \dint  u_1^+(x)\dint_{|x-y|\ge \eta} \frac {1}{|x-y|^{N+1-2s}} \,dy\,dx  \leq   C \lambda  \dint  u_1^+(x)\, dx .
 \end{array}
$$
In both cases we get $ {\cal O}(\lambda^2)$ for an appropriate $\nu$, which proves \eqref{Lipschitz}  as wanted.

\medskip

\noindent (iii) Estimate \eqref{Lipschitz} shows a property of one-sided Lipschitz continuity in time from above for the space integral of the LHS. On the other hand, we also have the inequality
$$
\dint \varphi_1 H(g) (t) dx \le \dint \varphi_1 (g-1)^2 (t) dx \le \dint u^+_1(x)^2 dx \leq C\lambda^2, \quad \forall t.
$$
Inserting this into the energy inequality \eqref{Lipschitz} we arrive at
$$
\dint_{T_1}^{T_2} {\cal B}(u_1^+,u_1^-)dt \le C \lambda^2, \quad {\rm for}
\quad -4<T_1<T_2<0.
$$
We are now in the position to apply Steps 2 and 3 of the proof of Section 4 of \cite{CChV}   with only technical changes that we will omit, and we  will thus get Lemma  \ref{lem.3layer}. Note that the quadratic (or at least superlinear) dependence on $\lambda$ in Formula \eqref{Lipschitz}  is absolutely necessary for the proof to work.
\qed


\section{Oscillation reduction result}

We are now ready to prove the strong oscillation reduction result, Lemma \ref{reg.2}.
Though the general idea of the argument follows closely Section 5 of \cite{CChV},  an interesting modification is needed to accomodate the nonlinearity of the equation: in doing the scaling of the solutions in each iteration we will now find solutions of a family of related equations, and we have to use the modified estimates of Section \ref{sec-mod}, that hold for that family. At some step we will have to apply the results of the previous section to such a family. This is an easy further verification that we will leave to the reader.

Instead of Lemma \ref{reg.2}, we will prove the following version that will be the one used later in the regularity argument. Let us fix some notation. For $\lambda$ as in the previous section, we define for any $\ve>0$
$$
\psi_{\ve,\lambda}(x)=((|x|-1/\la^{4/s})^\ve-1)_+ \quad \mbox{if \ } \quad |x|\ge \lambda^{-4/s},
$$
and zero otherwise.

\begin{lem}\label{lem10.1} Given $\rho>0$ there exist $\ve>0$ and $\mu_1$ such that for any solution of the FPME in $(-4,0)\times \ren$ satisfying
\begin{equation}
0\le u \le 1+\psi_{\ve,\lambda},
\end{equation}
and assuming that
\begin{equation}
|\{u<\vp_0\}\cap (-4,-2)\times B_1|>\rho\,,
\end{equation}
then we have
\begin{equation}
\sup_{(-1,0)\times B_1} u \le 1-\mu_1.
\end{equation}
\end{lem}

 Note: In the next section we will take $\rho$ equal to $\delta_0$ as defined after the statement of  Lemma 5.2.

\noindent {\sl Proof.} (i) Consider $j_0=|(-3,0)\times B_3|/\gamma,$ with the value of $\gamma$  given in Lemma \ref{lem.3layer}. We fix $\ve>0$ small enough so that
$$
\frac{(|x|^\ve-1)_+}{\lambda^{2j_0}}\le (|x|^{s/4}-1)_+
$$
for all $x$. We may take $\ve=(s/4)\lambda^{2j_0}$ for instance. For $j\le j_0$ we consider
the sequence of functions defined iteratively by
\begin{equation}
u_{j+1}(x,t)=\frac1{\lambda^2}(u_j(x,t)-(1-\lambda^2))\,,
\end{equation}
starting from $u_0(x,t)=u(x,t)$, the solution of the FPME under consideration. By induction we assume that
$$
(u_j)_+(x,t)\le 1+ \frac1{\lambda^{2j}} \psi_{\ve,\la}(x)\,, \quad (x,t)\le (-3,0)\times \ren.
$$
So for $j\le j_0$ we have $u_j\le 1+\psi_\la$. We recall that $\lambda$ is fixed and small.

\medskip

\noindent (ii) Next, we have to check the equation satisfied by
$u_j$, that turns out to be
\begin{equation}
u_t=\nabla (D_j(u)\nabla {\cal L}(u))
\end{equation}
with a diffusion coefficient $D_j$ defined inductively by the rule
\begin{equation}
D_{j+1}(s)=D_{j}(\la^2 s + (1-\la^2)),
\end{equation}
so that  $D_j(u_j)=D_{j+1}(u_{j+1})$ holds. These are the type of diffusion coefficients for which
we have proved the modified estimates of Section \ref{sec-mod}.  We have to make the observation that, with the notations of that section, we have  $d_2u_{j+1}=\lambda^2 u_{j+1}= u_j-(1-\lambda^2)$ so that the integrals in $\cal Q$ do not change their form, in particular the formulas of the type \eqref{delicate}, that depend on the integral
$$
\int \nabla L(x-y)(u_{j+1}(y)-\vp_k(y))^-\,dy
$$
are estimated by the same constants. Note that  $\int_{\ren} \nabla L(x-y)\,C\,dy=0$ by antisymmetry of $\nabla L$.

(iii) By construction, the measure $|\{u_j\le \vp_0\}\cap (-4,0)\times B_1|$ is increasing; hence, it is bigger than $\rho$ for every $j$. We can apply Lemma \ref{lem.3layer} inductively to $u_j$. As long as \\ $|\{u_k>\vp_2\}\cap ((-2,0)\times \ren)|\ge \delta$ for $k=1,2,\cdots,  j+1$, we have
$$
|\{\vp_0\le u_{j+1}\le \vp_2\}\cap(-3,0)\times \mathbb R^N\}|\ge \gamma.
$$
Therefore,
$$
\begin{array}{l}
|\{u_{j+1}>\vp_2\}| \le |\{u_{j+1}>\vp_0\}|-\gamma
\le |\{u_{j}>\vp_2\}|-\gamma\le |(-3,0)\times B_3|-j\gamma\,.
\end{array}
$$
This cannot be true up to $j_0$, hence there must exist $j\le j_0$ such that
$$
|\{u_{j}>\vp_2\}\cap ((-2,0)\times\ren) | < \delta\,.
$$
We can then apply the first lemma to the next step, $u_{j+1}$. Indeed,
$$
u_{j+1} \le 1+\psi_\lambda\le 1+\psi_1 \quad \mbox{on } \ (-3,0)\times \ren
$$
and
$$
\begin{array}{l}
|\{u_{j+1}>1/2\}\cap ((-2,0)\times B_2)| \le |\{u_{j+1}>\vp_0\}\cap ((-2,0)\times B_2)|\\[5mm]
\le |\{u_{j}>\vp_2\}\cap ((-2,0)\times \ren)|< \delta\,,
\end{array}
$$
and then Lemma \ref{reg.1} implies that
$$
u_{j+1} \le 1-\mu \quad \mbox{ on \ } \quad (-1,0)\times B_1.
$$
This gives the result with $\mu_1=\lambda^{2j_0} \mu .$ \qed

\


\section{End of proof of regularity for $s<1/2$}\label{sec.finalCalpha}

The whole technical machinery is in place, and we are ready to prove H\"older regularity by means of the iterative process outlined in Section \ref{sec.holder}. We take any point $P_0=(x_0, t_0) \in \ren\times(0,\infty)$ and prove that $u$ is $C^\alpha$ around $P_0$
with an exponent that depends only on the parameters $N$ and $s$ of the equation, and a H\"older constant that depends also on the $L^\infty$ norm of $u$ and a lower bound of $t_0$.

We start with some reductions. There is no loss of generality in assuming that $u$ is bounded in the cylinder $Q$ since we know by Theorem \ref{th:L-inf} that $u$ is bounded in any strip of the form $\ren\times (\tau,\infty)$ with $\tau>0$. Moreover, by scaling we may assume that $t_0>4$. It will be then convenient to make a space-time translation and put $P_0=(0,0)$ assuming that the domain of definition of $u$ contains the strip $S_4=\ren\times [-4,0]$. By scaling we may assume also that $0\le u(x,t)\le 1$ in $S_4$.

\noindent Consider now a positive constant $K < 1/4$ such that the growth of the wings is controlled as follows:
\begin{equation}\label{wing.growth}
 \frac1{1-(\mu_1/2)}\psi_{\lambda,\ve}(Kx)\le \psi_{\lambda,\ve}(x).
\end{equation}
 The coefficient $K$ depends only on $\lambda, \mu_1$ and $\ve>0$.
 The parameters are as in the last section. The iteration that we will perform offers two possibilities.

 \noindent $\bullet$ {\sl Alternative 1. Regularity at a degenerate point.} Suppose that we can apply Lemma \ref{lem10.1} repeatedly  because the lowering of the oscillation may be assumed to happen always from above. We consider then the sequence of  functions  defined in the strip  $S_3=\ren\times (-4,0)$ by
  \begin{equation}\label{reg.scal}
 u_{j+1}(x,t)=\frac1{1-(\mu_1/4)}\, u_j(Kx, K_1 t), \quad K_1=\frac{K^{2-2s}}{1-(\mu_1/4)}\,.
 \end{equation}
 Note that this time the $u_j$'s are all of them solutions of the same equation.  According to the running assumption, and using \eqref{wing.growth}, we can apply Lemma \ref{lem10.1} at every step so that we have $u_{j}(x,t)\le 1-\mu_1$ in the cylinder $Q_1=B_1\times (-1,0)$ for every $j\ge 1$. In view of the scaling \eqref{reg.scal}, this implies H\"older regularity around the point  $(0,0)$, where the solution necessarily takes the degenerate value $u=0$ in a continuous way.

 \medskip

  \noindent $\bullet$ {\sl Alternative  2. Regularity at points of positivity.} After some steps of the iteration the assumption on the measure of the set $\{u_j>1/2\}$ made in Lemma \ref{lem10.1}  fails. Then, we are in the situation where the oscillation is reduced from below thanks to Lemma \ref{reg.1b}, which pulls the solution uniformly up from zero in a smaller cylinder. Then the equation is no longer degenerate, because after that
  step we have
  $$
  0<\mu'\le u_j(x,t)\le 1,
  $$
  in the cylinder $B_1\times (-1,0)$. Scaling the situation we will be in the conditions of the equation with
  diffusivity $D(u)$ mentioned above,  to which we apply either the modification of Lemma \ref{lem10.1} or the modification of Lemma \ref{reg.1b}. In fact, we can apply the modification of  Lemma \ref{lem10.1} both from above and from below since the degeneracy has disappeared. In this way, we obtain H\"older regularity at a point $P_0$ where $u(P_0)>0$.


\section{$C^\alpha$ regularity for $s> 1/2$}\label{s-large}

Here we want to prove the regularity result of Theorem \ref{mainthm} for $1/2<s<2$.
 If we  try to re-do the proof of the technical oscillation lemmas listed in Section \ref{sec.holder}, we find a convergence problem that is already apparent in Lemma \ref{reg.1}. Indeed, the bulk of the proof contained in Section \ref{sec51} works without modification, and  a problem was found only in the last estimate of Subsection \ref{sec51}.4, regarding integral ${\cal Q}_3$ in an outer region, where the decay of the factors that we have been examining does not provide for a uniform convergence of the integrals if $s\ge 1/2$.

\subsection{Analysis of the difficulty}

A possible solution is to make use of the known fact  that $u(x,t)$ is an  $L^1$ function in $x$
uniformly in $t$, in order to bound the integral of the $y$ terms in ${\cal Q}_3$, with integrand $\nabla L(x-y)u_k^-(x) $, since $L(x-y)$ is bounded for large $|y|$. Indeed, we can put together the different parts of  ${\cal Q}$ in the outer domain, and after dropping the contribution of the term in $u(x)$ since it is zero by the antisymmetry property of $\nabla L$, we get
$$
\begin{array}{l}
\left|\diint u_k^+(x)\frac{\nabla \vp(x)}{\vp(x)} \nabla L (x-y)\,u(y)\,dxdy\right|\\
\le C \|u(\cdot,t)\|_1 \left(\dint u_k^+(x)dx\right)\le C^{k}\|u(\cdot,t)\|_1 \left(\dint (u_{k-1}^+(x))^2dx\right).
\end{array}
$$
This solves the problem in one application of Lemma \ref{reg.1}, but then $\delta$ and $\mu$ would depend on $\|u(t)\|_1$. However, to obtain the $C^\alpha$ regularity result we have seen in the preceding Section \ref{sec.finalCalpha} that we need to iterate this and the other oscillation  lemmas, we want to rescale and repeat, and then a difficulty re-appears, because we will keep expanding the $u$ and the $x$ and therefore expanding
the integral at every step, so the constants will be ruined in the iteration. We need a way
to control such  behaviour.

 It will be convenient to examine the whole part of ${\cal Q}(u^+,u)$ that contains the difficulty, i.e.,
\begin{equation}
   \diint u_k^+(x)\frac{\nabla \vp_k(x)}{\vp_k(x)} \nabla L (x,y) \,u(y)\,dxdy\,.
\end{equation}
 As we dilate and repeat in the iteration scheme, the term $ \nabla L(x-y) u(y)$
 also starts to build up as $y$ tends to infinity. On the other hand, the integrability in $y$ at infinity is lost if $s\ge 1/2$, since in that case $\nabla L$ decays  like
$$
 |\nabla L| \sim |y|^{-(N-2s+1)}
 $$
and this is not good enough. However, the good news is that
 $$
 \nabla \cdot(\nabla L) \sim |y|^{-(N+2-2s)}\,,
 $$
(valid for all second derivatives) which is integrable as $|y|\to\infty$. Noting that $\nabla L$ is integrable for $y\sim 0$ if $s>1/2$, we conclude that
$$
V(x,t): = (\nabla L (x,y,t)) \ast_y u(y,t)
$$
has bounded H\"older seminorm. Therefore, it would be enough to control it at just a point, for instance at $x=0$ (for an interval of times).

\medskip

\subsection{ The transport approach for the case $s>1/2$}

The technical way to make use of the last observation is to perform by a change of coordinates $x'=x-S(t) $ that introduces a transport term to counter the difficult term $ \int \nabla L(y)\,\Psi(y)\, u(y)\,dy. $ To be precise, we define
\begin{equation}\label{def.vel}
S(t)=\dint_0^t \vec{v}(t)\,dt, \qquad
\vec{v}(t)=\dint \nabla L\, (y)\,u(y,t)\,dy\,,
\end{equation}
and we observe that $|\vec{v}(t)|$ depends also on $u$, and that  $|\vec{v}(t)|\le C<\infty $ since $u(y,t)$ is in $L^1_y$ uniformly in $t$. Indeed, the value of $\vec{v}(t)=\vec{v}(t; u)$ is only controlled by the integral of $u$ is space (i.\,e., the mass of $u(t)$).
Next, we introduce change of variables
$$
(x,t) \mapsto (x',t'):= (x-S(t),t)\,,
$$
and we write the equation for $u$ with respect to the new variables, $u(x,t)=\tilde u(x',t')$. The RHS does not change since  we are performing a space translation for fixed time. However, the time derivative in the LHS transforms as follows:
$$
u_t(x,t) = u_{t'}(x'+S (t'),t')+(\nabla u)S'(t')= \tilde u_{t'} + \vec{v}\cdot\nabla_{x'}   \tilde u.
$$
The new term is what we are aiming at.  The equation takes the convection-diffusion form
\begin{equation}\label{eq.diff.conv}
\tilde u_t+\vec{v}\cdot\nabla_{x'}   \tilde u = \nabla (\tilde  u\,\nabla {\cal L} \tilde u)\,.
\end{equation}
In the sequel we will write  $t$ for $t'$ and $u$ instead of $\tilde u$ without fear of confusion. The space variable is still written $x'$. Next, we pass the term $\vec{v}\cdot\nabla u$ to the RHS and  multiply by $\log((u/\vp)\vee 1)$, as we did in Section \ref{sec51}, to obtain  the energy formula. We observe that in this case the RHS contains an extra term of the form
$$
I= \diint  \nabla \log( (u/\vp)\vee 1) \,\vec{v}(t) u(x',t)\,dx'dt\,.
$$
This integral must be computed only in the region where $u>\vp$, and in that case $(u/\vp)\vee 1=u/\vp=1+(u_k^+/\varphi_k)$, so that
$$
I= \dint dt\dint_{u>\vp_k}  \nabla u_k^+  \, \vec{v}(t) dx' -
\dint dt\dint_{u>\vp_k} u_k^+ \frac{\nabla \vp }{\vp}  \, \vec{v}(t) dx' =I_1-I_2
$$
The first integral vanishes, and the second is precisely the troublesome  term:
$$
I_2=
\dint dt \dint_{u>\vp_k} dx (u_k^+(x')/\varphi_k(x')) \nabla \varphi_k(x')
 \dint \nabla L(y') u(y',t)\,dy'\,,
$$
that disappears in this way from the calculation.  Alternatively, the disappearance of the bad term in the energy calculation in the new variables can be easily seen if we write the equation for $\tilde u(x',t)$  in the more symmetrical form
\begin{equation}\label{trans.eq}
\tilde u_t=\nabla_{x'} \left(\tilde u  \left(\int \{\nabla L(y-x')-\nabla L(y)\}(\tilde u(y)-\tilde u(x'))\,dy\right)\right)\,,
\end{equation}
to be interpreted in the same weak form, or weak energy form, that we used for $u(x,t)$. In any case, this allows  to prove Lemma \ref{reg.1} also for $s\in (1/2,1)$ if we work in the new coordinates, and the constants involved in the result do not depend on the $L^1$ norm of the solution. The price to pay is that the slope of the distorted space variables does depend on the $u$-integral. So, in the first step of the iteration process we have shown how to transfer the difficulty from a numerical term to a geometrical distortion.

In order to sum up the result, let us introduce the bound $M=1\vee \sup_{t>0} \vec{v}(t)$,
that depends only on $u$ via the norm $\sup_t\|u(\cdot,t)\|_1$.

\begin{lem} \label{reg.1Mod} Let $1/2< s<1$ and let $u$ be a solution of the FPME under the assumptions of Lemma {\rm \ref{reg.1}}. Let us perform the above change of variables so that $\tilde u(x',t')$ is defined in
smaller cylinder $Q_L$ where $L=4/(M+1)$. Then the result of Lemma {\rm \ref{reg.1}} is true for $\tilde u$, with conclusion holding in a smaller cylinder $Q_{1/M}$;  $\delta$ may depend also on $M$.
\end{lem}

Thanks to \ref{trans.eq},  it is then immediate to see that the modified Lemma \ref{reg.1.mod}, as well as the pull-up Lemma \ref{reg.2} are also true if stated in the form that we have used for Lemma \ref{reg.1Mod}. A bit more of attention to the details will show that the stronger reduction Lemma \ref{lem10.1} also holds, since the iterations do not change the scaling in space and time.

\subsection{Analysis of the transport term in the final iteration}

When we try to perform again the iteration procedure of Section \ref{sec.finalCalpha}, one of the alternatives is repeated scaling around a degenerate point. In that case the iterations take the form
\begin{equation}
 u_{j+1}(x,t)=\frac1{1-(\lambda_*/4)}\, u_j(Kx, K_1 t), \quad K_1=\frac{K^{2-2s}}{1-(\lambda_*/4)}\,.
 \end{equation}
 that we may sum up as
$$
u_1(x_1,t_1)=A\,u(x,t) , \quad x_1= B x, \quad t_1=C t.
$$
where $A<1$, $B<1$ and $C=B^{2-2s}A$, so that the same equation will be satisfied after the
change of scale. We propose here to do the same iteration for the solution $\tilde u$ in terms of the variables $x'$ and $t$. The equation will then take the modified form \eqref{trans.eq}, that will be satisfied again by the iterates, just as it is written.  It is true that the velocity $\vec{v}(t)$ will change from iteration
to iteration according to the rule
$$
\vec{v_1}(t)= \frac{C}{B}\,\vec{v}(Ct)=B^{1-2s} \vec{v}(Ct),
$$
which follows both from the geometric transformation, and from the definition of $\vec{v}$ in \eqref{def.vel}.
Therefore, after the first geometrical transformation such repeated iterations  conserve the same correspondence for all subsequent steps. In other words, the geometrical transformation done in the first step will hold for all remaining steps: if the set of coordinates at that moment is  $(x_n,t_n)$, we obtain a set of newly distorted  coordinates $(x'_{n},t_n)$ by the formula
$$
x'_{n}(t)=x_{n}(t)-S_n(t), \qquad S_v'(t)=\vec{v}_{n}(t)
$$
then this is just a scaled version of the original transformation for $n=0$. Summing up, since the contractions
in the upper bound for $u$ happen with a constant rate $1-\mu$ in cylinders that shrink in space and time also  with a fixed rate, we conclude in a standard way $C^\alpha$-regularity with respect to the transformed variables. But since the coordinate  transformation is done only once and is Lipschitz continuous, this means the same type of H\"older regularity for $u$ with respect to the original coordinates $(x,t)$. Of course, the Lipschitz constant of the transformation depends on $M_1=\sup_t \|u(t)\|_{L^1_x}$.

\medskip

The analysis of the second alternative is easier since we are converging along the iterations to an equation with constant diffusivity coefficient. We leave the easy details to the reader. This ends the proof of Theorem
\ref{mainthm} for $1/2<s<2$. \qed

\section{Extension of the existence theory}\label{sec.ext}

After these results, we can extend the existence theory to all nonnegative and integrable initial data.

\begin{thm} For every $u_0\in L^1(\ren)$, $u_0\ge 0 $, there exists a continuous weak solution of the FPME in the following sense: there exists a function $u(x,t)$, continuous and nonnegative in $Q=\RR^n\times (0,T)$ such that \
$$
u\in L^\infty(0,\infty:L^1(\ren)\cap L^\infty(\ren\times (\tau,\infty) \ \mbox{ for all \ } \tau>0\,,
$$
$$
{\cal K}(u)\in L^1(0,  T: W^{1,1}_{loc}(\RR^n)), \qquad u\,\nabla{\cal K}(u)\in L^1(Q_T)
$$
and the identity
\begin{equation}\label{identity}
\iint u\,(\eta_t-\nabla {\cal K}(u)\cdot\nabla\eta)\,dxdt+ \int
u_0(x)\,\eta(x,0)\,dx=0
\end{equation}
holds for all  continuously differentiable test functions $\eta$ in $Q_T$  that are compactly supported in the space variable and vanish near $t=T$.
\end{thm}

\noindent {\sl Proof.} (i) We take a sequence of initial data $u_{0n}$ that are nonnegative, smooth and decaying at infinity as required in \cite{CaVa09}. We assume that $\|u_{on}\|_1\le C_1$ for all $n\ge 1$ and $u_{0n}\to u_0$ in $L^1(\ren)$. Then there exist solutions $u_n(x,t)$ and there are estimates like $u_n\ge 0$ and
$$
\int u_n(x,t)\,dx=\|u_{0n}\|_1\le C_1  \qquad \forall n\ge 1, \ t>0\,.
$$

(ii) We also know from Theorem \ref{th:L-inf} that the solutions are bounded for $t\ge\tau>0$ with bound that depends only on the $L^1$-norm of the initial data
$$
0\le u_n(x,t)\le C(n,s)\|u_{0n}\|_1^{\gamma}\,t^{-\alpha}\le C\,C_1^{\gamma}\,t^{-\alpha}.
$$
We also know that uniformly bounded nonnegative solutions are $C^\beta$ smooth for some  $\beta>0$ with
uniform H\"older constants, so that after passing to a subsequence we have
$$
u_n(x,t)\to u(x,t) \quad \mbox{uniformly on compact subsets of $Q$}\,.
$$

(iii) We now observe that for every $t>0$ we have $u_n(t)\in L^1\cap L^\infty$ which together with the
Sobolev embedding (or Riesz embedding) gives
$$
H(u_n(t))\in L^r(\ren) \qquad  \forall \ r\ge n/(n-s)
$$
and the embedding is compact into $L^r_{loc}(\ren)$.
Indeed, if $u_n\in L*p$ then
$$
 \frac 1{r}=\frac 1{p}-\frac{s}n, \qquad p\ge 1.
$$
In practice we will need the estimate
$$
\int H(u_n(t))^2\dx\le C \|u_n(t)\|_p^2, \quad \frac 1{2}=\frac 1{p}-\frac{s}n\,.
$$
But since
$$
\|u_n(t)\|_p^p\le  \|u_n(t)\|_1 \|u_n(t)\|_\infty^{p-1} \le
 C\|u(t)\|_1^{1+(p-1)\gamma}t^{-\alpha (p-1)}\,,
$$
we get the decay estimate for this energy in the form
$$
\int H(u_n(t))^2\dx\le C\,\|u(t)\|_1^{2(1+(p-1)\gamma)/p}\,t^{-2\alpha (p-1)/p}=
 C\,\|u(t)\|_1^{\sigma} t^{\lambda-1},
$$
with
$$
\lambda=1-\frac{2\alpha (p-1)}{p}=1-\frac{2n}{(n+2-2s)}\frac{(n-2s)}{2n}=1-\frac{n-2s}{n+2-2s}=\frac{2}{n+2-2s}.
$$

(iv)  We also have for every $t_2\ge t_1\ge 0$ the energy inequality
$$
\frac12 \int H(u_n(t_2))^2\,dx + \int_{t_1}^{t_2}\int u_n\,|\nabla K u_n|^2\,dxdt\le \frac12\int H(u_n(t_1))^2\,dx
$$
Therefore, we have a uniform estimate for the integral $\iint u_n\,|\nabla K u_n|^2\,dxdt$ from $t=\tau>0$ up to $t=$ infinity, and we get the same decay rate.

(v) Next, we address the step of passing to the limit of the weak formulation of the solutions $u_n$ when the test function
$\eta\in C^1(Q)$ that  compactly supported in $Q$ away from $t=0$:
\begin{equation}
\iint u_n\,(\eta_t-\nabla {\cal K}(u_n)\cdot\nabla\eta)\,dxd+ \int
u_{0n}(x)\,\eta(x,0)\,dxt=0.
\end{equation}
There is no problem in the convergence of $\iint u_n\,\eta_t\,dxdt$. If $\eta$ is supported away from $t=0$
then the second term converges since we can use the
compactness of the map $u\mapsto \nabla K(u)$ from $L^1\cap L^\infty $ into $L^r$ for some $r$

(vi) We now want to check that the initial data are taken.  For that we show first that
$u\nabla K u\in L^1_{x,t}$ also near $t=0$. Here is the calculation in a cylinder $Q_k=\ren\times (t_k,t_{k-1})$ with $t_k=2^{-k}$:
$$
\begin{array}{l}
\diint_{Q_k}u\,|\nabla K(u)|\,dxdt
\le  (\diint_{Q_k}  u\,dxdt)^{1/2}(\diint_{Q_k} u\,|\nabla K(u)|^2\,dxdt)^{1/2}\\[14pt]
\le  (\|u_0\|_1t_{k-1})^{1/2}(\frac12\dint H(u(t_k))^2\,dx)^{1/2}\le Ct_k^{1/2} t_k^{(\lambda-1)/2}=C\,t_k^{\lambda/2}.
\end{array}
$$
With this it is possible to pass to the limit in the term $\iint u_n\,\nabla {\cal K}(u_n)\cdot\nabla\eta\,dxdt $ for $\eta$ that does not vanish near $t=0$ by estimating the integral for small $t$ as uniformly small, and then proving the convergence for $t\ge \tau>0$ using the known regularity. This proves that the definition of solution according to Formula \eqref{identity} is true.

(vii) We can also check that the initial data are taken as traces.
This depends on a uniform estimate of the difference $u(\tau)-u_0$ for small $\tau$ in some adequate norm.
We have that for a $C^1_0$ test function $\eta(x)$
$$
\begin{array}{l}
\dint |u(\tau)-u_0|\,\eta\,dx\le \dint_0^\tau\int |\partial_t u |\eta \,dxdt\le \dint u\,|\nabla K(u)|\,|\nabla \eta|\,dxdt \\
\le C (\diint_K  u|\nabla \eta|^2\,dxdt)^{1/2}(\diint_K u\,|\nabla K(u)|^2\,dxdt)^{1/2}\le  C\,\tau^{\lambda/2}.
\end{array}
$$
This is what makes $u(t)\to u_0$ in $W^{-1,1}_{loc}$ in the limit.  Or we may use J. Simon's compactness results \cite{Simon},  since we are actually proving that $u_t\in L^1_t(W^{-1,1}_x)$.\qed

\


\noindent {\large\bf Acknowledgments}.   Part of the present work was done during stays of F. Soria and J. L. V\'azquez at Univ. of Texas in Austin, and of \ L. Caffarelli and J. L. V\'azquez at MSRI, Berkeley, CA. L. Caffarelli has been funded by NSF Grant DMS-0654267  (Analytical and Geometrical Problems in Non Linear Partial Differential Equations), F. Soria by Spanish Grant MTM2010-18128, and J. L. V\'azquez by Spanish Grant MTM2008-06326-C02-01.

\vskip 1cm


\bibliographystyle{amsplain} 

\

{\sc Addresses:}

\medskip

{\sc Luis A. Caffarelli}\newline
School of Mathematics, Univ. of Texas at Austin,
1 University Station, C1200, Austin, Texas 78712-1082. \newline
Second affiliation: Institute for Computational Engineering and Sciences.\newline
e-mail: caffarel@math.utexas.edu

\medskip

{\sc Fernando Soria}\newline
Departamento de Matem\'{a}ticas, Universidad Aut\'{o}noma de Madrid, 28049
Madrid, Spain.  Second affiliation: Institute ICMAT. \newline
e-mail: fernando.soria@uam.es

\medskip

{\sc Juan Luis V{\'a}zquez}\newline
Departamento de Matem\'{a}ticas, Universidad Aut\'{o}noma de Madrid, 28049
Madrid, Spain.  Second affiliation: Institute ICMAT. \newline
e-mail: juanluis.vazquez@uam.es

\vskip 1cm

2000 {\bf Mathematics Subject Classification.} 35K55, 35K65, 76S05.

{\bf Keywords and phases.} Porous medium equation, fractional Laplacian, nonlocal operator, regularity.

\vskip1cm
\end{document}